\documentclass[hidelinks,onefignum,onetabnum]{siamart250211}

\usepackage{lipsum}
\usepackage{amsfonts}
\usepackage{graphicx}
\usepackage{epstopdf}
\usepackage{algorithmic}
\ifpdf
  \DeclareGraphicsExtensions{.eps,.pdf,.png,.jpg}
\else
  \DeclareGraphicsExtensions{.eps}
\fi

\newsiamremark{remark}{Remark}
\newsiamremark{hypothesis}{Hypothesis}
\crefname{hypothesis}{Hypothesis}{Hypotheses}
\newsiamthm{claim}{Claim}

\headers{Parallelisation of PDEs via representation theory}{S. Olver}

\title{Parallelisation of  partial differential equations via representation theory
\thanks{Submitted to the editors 15 October 2025.
\funding{This work was supported in part by the EPSRC grant EP/T022132/1 “Spectral element methods for fractional differential equations, with applications in applied analysis and medical imaging”.}}}

\author{Sheehan Olver\thanks{Imperial College, London, UK 
  (\email{s.olver@imperial.ac.uk}, \url{https://www.ma.imperial.ac.uk/\string~solver/}).}}

\usepackage{amsopn}
\DeclareMathOperator{\diag}{diag}

\usepackage{pgf,ytableau}

\usepackage{xcolor}
\usepackage{arydshln}
\usepackage{mathdots}

\usepackage{todonotes}

\def\addtab#1={#1\;&=}

\def\meeq#1{\def\ccr{\\\addtab}
 \begin{align*}
 \addtab#1
 \end{align*}
  }  
  
  \def\leqaddtab#1\leq{#1\;&\leq}

\def\vc#1{\mbox{\boldmath$#1$\unboldmath}}

\def\pr(#1){\left({#1}\right)}
\def\br[#1]{\left[{#1}\right]}
\def\fbr[#1]{\!\left[{#1}\right]}
\def\set#1{\left\{{#1}\right\}}
\def\ip<#1>{\left\langle{#1}\right\rangle}
\def\iip<#1>{\left\langle\!\langle{#1}\right\rangle\!\rangle}

\def\fpr(#1){\!\pr({#1})}

\def\floor#1{\left\lfloor#1\right\rfloor}

\def\mapengine#1,#2.{\mapfunction{#1}\ifx\void#2\else\mapengine #2.\fi }

\def\map[#1]{\mapengine #1,\void.}

\def\mapenginesep_#1#2,#3.{\mapfunction{#2}\ifx\void#3\else#1\mapengine #3.\fi }

\def\mapsep_#1[#2]{\mapenginesep_{#1}#2,\void.}

\def\vcbr{\br}

\def\bvect[#1,#2]{
{
\def\dots{\cdots}
\def\mapfunction##1{\ | \  ##1}
\begin{pmatrix}
		 \,#1\map[#2]\,
\end{pmatrix}
}
}

\def\vect[#1]{
{\def\dots{\ldots}
	\vcbr[{#1}]
}}

\def\vectt[#1]{
{\def\dots{\ldots}
	\vect[{#1}]^{\top}
}}

\def\Vectt[#1]{
{
\def\mapfunction##1{##1 \cr} 
\def\dots{\vdots}
	\begin{pmatrix}
		\map[#1]
	\end{pmatrix}
}}

\def\E{{\rm e}}

\def\D{{\rm d}}

\def\tF_#1{{\tt F}_{#1}}

\def\tFC_#1{{\tt T}_{#1}}

\def\secref#1{Section~\ref{Section:#1}}

\def\lmref#1{Lemma~\ref{Lemma:#1}}

\def\thref#1{Theorem~\ref{Theorem:#1}}

\def\qand{\quad\hbox{and}\quad}
\def\qqand{\qquad\hbox{and}\qquad}

\def\Abstract#1\par{\begin{abstract}#1\end{abstract}}
\def\Keywords#1\par{\begin{keywords}{#1}\end{keywords}}
\def\Section#1#2.{\section{#2}\label{Section:#1} }
\def\Appendix#1#2.{\appendix \section{#2}\label{Section:#1} }

\def\Subsectionl#1#2.{\subsection{#2}\label{subsec:#1}}
\def\Subsection#1.{\subsection{#1}}

\def\Subsubsection#1.{\subsubsection{#1}}

\def\Problem#1#2\par{\begin{problem}\label{Problem:#1} #2\end{problem}}
\def\Theorem#1#2\par{\begin{theorem}\label{Theorem:#1} #2\end{theorem}}
\def\Conjecture#1#2\par{\begin{conjecture}\label{Conjecture:#1} #2\end{conjecture}}
\def\Proposition#1#2\par{\begin{proposition}\label{Proposition:#1} #2\end{proposition}}
\def\Definition#1#2\par{\begin{definition}\label{Definition:#1} #2\end{definition}}
\def\Corollary#1#2\par{\begin{corollary}\label{Corollary:#1} #2\end{corollary}}
\def\Lemma#1#2\par{\begin{lemma}\label{Lemma:#1} #2\end{lemma}}
\def\Example#1#2\par{\begin{example}\label{Example:#1} #2\end{example}}
\def\Remark #1\par{\begin{remark*}#1\end{remark*}}

\def\figref#1{Figure~\ref{fig:#1}}

\def\Figurew[#1]#2#3\par{
\begin{figure}[tb]
\begin{center}{
\includegraphics[width=#2]{Figures/#1}}
\end{center}
\caption{#3}\label{fig:#1} 
\end{figure}
}

\def\Figure[#1]#2\par{
\begin{figure}[tb]
\begin{center}{
\includegraphics{Figures/#1}}
\end{center}
\caption{#2}\label{fig:#1} 
\end{figure}
}

\def\Figurefixed[#1]#2\par{
\Figurew[#1]{0.48 \hsize}{#2}\par
}

\def\Figuretwow#1#2#3#4\par{
\begin{figure}[tb]
\begin{center}{
\includegraphics[width=#3]{Figures/#1}\includegraphics[width=#3]{Figures/#2}}
\end{center}
\caption{#4}\label{fig:#1} 
\end{figure}
}

\def\Figuretwowframed#1#2#3#4\par{
\begin{figure}[tb]
\begin{center}{
\fbox{\includegraphics[width=#3]{Figures/#1}}\fbox{\includegraphics[width=#3]{Figures/#2}}}
\end{center}
\caption{#4}\label{fig:#1} 
\end{figure}
}

\def\Figuretwo[#1,#2]#3\par{
	\Figuretwow{#1}{#2}{0.48 \hsize}
		#3\par	
}

\def\Figuretwoframed[#1,#2]#3\par{
	\Figuretwowframed{#1}{#2}{0.48 \hsize}
		#3\par	
}

\def\Figurethreew#1#2#3#4#5\par{
\begin{figure}[tb]
\begin{center}{
\includegraphics[width=#4]{Figures/#1} \includegraphics[width=#4]{Figures/#2} \includegraphics[width=#4]{Figures/#3}}
\end{center}
\caption{#5}\label{fig:#1} 
\end{figure}
}

\def\Figurethree#1#2#3#4\par{
	\Figurethreew{#1}{#2}{#3}{0.3 \hsize}
		{#4}\par	
}

\def\Figurematrixfour#1#2#3#4#5\par{
\begin{figure}[tb]
\begin{center}{
\vbox{\hbox{\includegraphics[width= 0.48 \hsize]{Figures/#1} \includegraphics[width= 0.48 \hsize]{Figures/#2}}\hbox{\includegraphics[width= 0.48 \hsize]{Figures/#3}\includegraphics[width= 0.48 \hsize]{Figures/#4}}}}
\end{center}
\caption{#5}\label{fig:#1} 
\end{figure}
}

\def\elllRpz_#1{\ell_{#1{\rm z}}^{(\lambda,R),p}}

\def\sopmatrix#1{\begin{pmatrix}#1\end{pmatrix}}

\def\Proof{\begin{proof}}
\def\mqed{\end{proof}}

\gdef\reffilename{\jobname}

\outer\def\ends{ 
\end{document}
}

\newtheorem{problem}{Problem}
\theoremstyle{definition}
\newtheorem{example}{Example}

\def\vect[#1]{
{\def\dots{\ldots}
	\pr({#1})
}}

\begin{document}





\def\GZn{\hbox{GZ}(n)}
\def\Spec{\hbox{Spec}}
\def\Cont{\hbox{\rm Cont}}
\def\cont{\hbox{\rm cont}}
\def\kron{\hbox{\rm kron}}
\def\vec{\hbox{\rm vec}\,}
\def\mspan{\hbox{\rm span}\,}
\def\colspan{\hbox{\rm colspan}\,}
\def\GL{\hbox{\rm GL}}
\def\mdiag{\hbox{\rm diag}}
\def\Y{{\mathcal Y}}
\def\algref#1{Algorithm~\ref{Algorithm:#1}}
\def\bfc{{\bf c}}
\def\bfx{{\vc x}}
\def\bfv{{\bf v}}
\def\bfzero{{\bf 0}}

\let\addition=\relax

\def\smallpmatrix#1{\left(\begin{smallmatrix} #1 \end{smallmatrix}\right)}

\def\probsref[#1,#2]{Problems \ref{Problem:#1} and \ref{Problem:#2}}

\def\bbR{{\mathbb R}}
\def\bbZ{{\mathbb Z}}
\def\bbC{{\mathbb C}}
\def\bbP{{\mathbb P}}
\def\bfP{{\mathbf P}}
\def\bfQ{{\mathbf Q}}
\def\bfp{{\mathbf p}}
\def\bfpt{{\mathbf p}\!{}^\top}
\def\bfq{{\mathbf q}}
\def\bfe{\vc e}
\def\bfq{\vc q}
\def\bfu{\vc u}
\def\bfk{\vc k}
\def\bfb{\vc b}
\def\bfd{\vc d}
\def\bfc{\vc c}
\def\bfone{\vc 1}

\def\rmC{{\rm C}}
\def\rms{{\rm s}}

\def\fee{f_{\rm ee}}
\def\feo{f_{\rm eo}}
\def\foo{f_{\rm oo}}
\def\foe{f_{\rm oe}}
\def\fre{f_{\rm re}}
\def\fro{f_{\rm ro}}
\def\ff#1{f_{{\rm f}#1}}
\def\fe{f_{\rm e}}
\def\fo{f_{\rm o}}
\def\ft{f_{\rm t}}
\def\fh{f_{\rm h}}
\def\fv{f_{\rm v}}
\def\fs{f_{\rm s}}
\def\sigmat{\sigma_{\rm t}}
\def\sigmats{\sigma_{\rm ref}}
\def\sigmast{\sigma_{\rm rot}}
\def\sigmas{\sigma_{\rm s}}
\def\rhoP{\rho_{\rm P}}
\def\sigmaf{\sigma_{\rm f}}
\def\sigmah{\sigma_{\rm h}}
\def\sigmav{\sigma_{\rm v}}
\def\sign{{\rm sign}\,}

\def\Orth{\hbox{\rm O}}
\def\Ord{{\mathcal O}}

\def\appref#1{Appendix~\ref{Section:#1}}
\def\tabref#1{Table~\ref{Table:#1}}

\def\rhos{\rho_{\rms}}

\setlength{\dashlinedash}{0.5pt} 
\setlength{\dashlinegap}{1.3pt}    

%
%

\maketitle

\begin{abstract}
\addition{
Incorporating symmetries into the numerical solution of differential equations has been a mainstay of research over the last 40 years, however, one  aspect is less known and under-utilised:  discretisations of partial differential equations that commute with symmetry actions (like rotations, reflections or permutations) can be decoupled into independent systems solvable in parallel by incorporating knowledge from representation theory. 
We introduce this beautiful subject  via a crash course in representation theory focussed on hands-on examples for the symmetry groups of the square and cube, and its utilisation in the construction of so-called {\it symmetry-adapted bases}.  Schur's lemma, which is not well-known in applied mathematics, plays a powerful role in proving sparsity of resulting discretisations and thereby  showing that partial differential equations do indeed decouple. Using  Schrödinger equations as a motivating example, we demonstrate that a symmetry-adapted basis leads to a significant increase in the number of independent linear systems. Counterintuitively, the effectiveness of this approach is in fact greater for partial differential equations with {\it less} symmetries,  for example a Schrödinger equation where the potential is only invariant under permutations, but not under rotations or reflections.  We also explore this phenomenon  as the dimension of the partial differential equation becomes large, hinting at the  potential for significant savings in high-dimensions.}

\end{abstract}

	
\Section{intro}  Introduction.

\addition{
Partial differential equations (PDEs) in applied mathematics  are often posed in symmetrical geometries (disks, balls, squares, or cubes) with variable coefficients that are invariant under  symmetry actions (rotations, reflections, or permutation).  As a motivating example, consider a Schr\"odinger equation in the unit cube $[-1,1]^3$: find $u(x,y,z)$  satisfying Neumann boundary conditions---the normal derivative vanishes at the boundary---such that
\[
\Delta u + \underbrace{a(x,y,z)}_{\hbox{Invariant}} u = \underbrace{f(x,y,z)}_{\hbox{Anisotropic}},
\]
where $\Delta = \partial_{x}^2 + \partial_{y}^2 + \partial_{z}^2$ is the Laplacian and the {\it potential} $a(x,y,z)$ is a function that is invariant under  rotations, reflections and permutations (e.g., $a(x,y,z) = x^2 + y^2 + z^2$) or under only negation and permutations (e.g., $a(x,y,z) = (x-y)^2 + (x-z)^2 + (y-z)^2$), with the latter case corresponding to interacting one-dimensional particles.  The right-hand side $f$ and the solution $u$ do not necessarily have any symmetry properties, that is, they may be {\it anisotropic} with respect to rotations, reflections or permutations.
}

\addition{
 This paper addresses a fundamental observation: incorporating symmetries of PDEs in the right way into their discretisation  will decouple the problem into independent linear systems that can be solved in parallel. The key idea is to build in knowledge about the basic ways in which symmetry can present itself,  which exactly corresponds to the building blocks of representation theory: {\it irreducible representations}. To understand this procedure we must have a basic understanding of representation theory, a topic that is not widely known in the applied mathematics community.
 }
 
%
%
%

\addition{
Let us motivate this procedure with a much simpler geometry: the unit interval $[-1,1]$, which is invariant under reflection. Reflection symmetry is naturally connected to a decomposition into even and odd parts:
\[
f(x)\ =\ \underbrace{\fe(x)}_{\hbox{Even}}\ +\ \underbrace{\fo(x)}_{\hbox{Odd}}
\]
where $\fe$ is invariant under reflection ($\fe(-x) = \fe(x)$), whilst $\fo$ is anti-invariant under reflection ($\fo(-x) = -\fo(x)$).  That is, the symmetry action of reflection ($x \mapsto -x$) now takes one of two forms:  either multiplication by $1$ (in the even case) or $-1$ (in the odd case).
%
A decomposition according to parity extends naturally to a geometry that has reflection symmetry in two-dimensions such as the unit square $[-1,1]^2$:
\[
f(x,y)\ = \underbrace{\fee(x,y)}_{\hbox{Even--Even}}\  + \ \underbrace{\foe(x,y)}_{\hbox{Odd--Even}}\ +\ \underbrace{\feo(x,y)}_{\hbox{Even--Odd}}\ +\ \underbrace{\foo(x,y)}_{\hbox{Odd--Odd}},
\]
where each component is either invariant or anti-invariant  under  horizontal and vertical reflections, see a depiction in \figref{reflectiondecomposition}. That is, the symmetry actions of horizontal  ($(x,y) \mapsto (-x,y)$) or vertical ($(x,y) \mapsto (x,-y)$) reflection become either multiplication by $1$ or $-1$, with each symmetry class exhibiting different behaviour.
}

\Figurew[reflectiondecomposition]{\hsize}
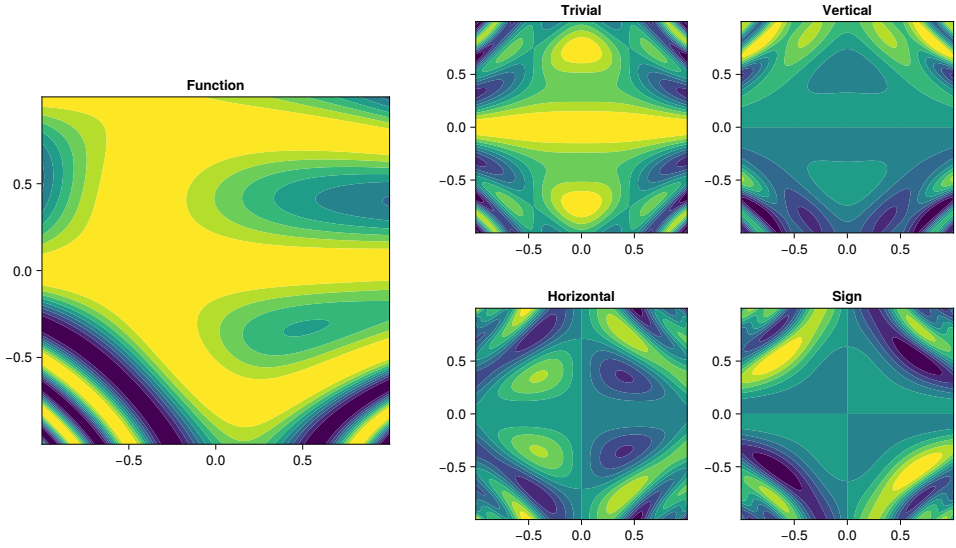
	A function on a square can be decomposed into four different symmetry classes according to their \addition{parity}.  On the left we plot an arbitrary function $f(x,y) = \cos(6y(y-1)(x-1/5)^2 + \sin(4y-1/10)\E^x)$ and on the right we plot the four components of its symmetric decomposition. Applying a \addition{horizontal or vertical} reflection leaves each term either invariant or swaps the sign. Certain PDEs can be parallelised across four independent solves according to this decomposition.

\addition{
The unit square has symmetries beyond horizontal and vertical reflection: adding in rotations it is invariant under the dihedral group $D_4$. This raises the question: can we break apart a function into different symmetry classes corresponding to  the different ways that $D_4$  can present itself?  We shall see that the corresponding decomposition is into 6 terms:
\[
f(x,y)\ =\ \underbrace{\ft(x,y)}_{\hbox{Trivial}}\ +\ \underbrace{\fre(x,y)}_{\hbox{Reflect}}\ +\ \underbrace{\fro(x,y)}_{\hbox{Rotate}}\ +\ \underbrace{\fs(x,y)}_{\hbox{Sign}}\ +\ \underbrace{\ff1(x,y) + \ff2(x,y)}_{\hbox{Faithful}},
\]
where the labels classify different symmetry classes, see a depiction in \figref{squaredecomposition}. The first four terms have {\it scalar} symmetries similar to the even--odd decomposition: $\ft(x,y)$ is invariant under reflections {\it and} rotations, $\fre(x,y)$  is anti-invariant under reflections but invariant under rotations, $\fro(x,y)$ is invariant under reflections and ant-invariant under $90^\circ$ rotations and $\fs(x,y)$ is anti-invariant under both   reflections and $90^\circ$ rotations.   That is to say rotations and reflections have become multiplication by $1$ or $-1$, with each symmetry class again exhibiting different behaviour.
But the last two terms ($\ff1$ and $\ff2$, corresponding to the ``Faithful representation'')  capture a phenomenon not present in a decomposition based on parity: they come as a pair, and applying a symmetry action like rotation or reflection will leave them inside a two-dimensional subspace, but this is not at all evident from their plots.
}

\Figurew[squaredecomposition]{\hsize}
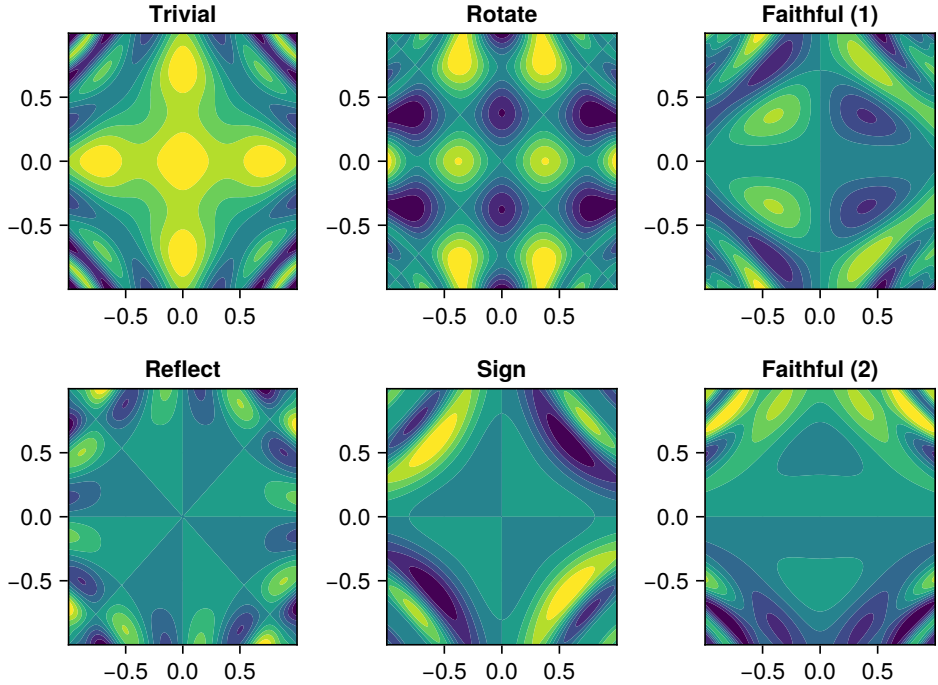
	A symmetric decomposition associated with the full symmetry group of the square (the dihedral group $D_4$) into six symmetry classes for the same function from \figref{reflectiondecomposition}.  The first two columns correspond to scalar irreducible representations where the functions are either invariant or flip sign whenever we apply \addition{a reflection or rotation}. The last column corresponds to a 2-dimensional irreducible representation,   \addition{where the symmetry is not at all visually obvious}.

\addition{
In order to make sense of the behaviour of reflections and rotations applied to $\ff1$ and $\ff2$ we  introduce the concept of a {\it symmetry-adapted basis}, which are bases that turn a symmetry action into multiplication.  The simplest example of a symmetry-adapted basis are in fact the monomials, which capture reflection symmetry: assuming the series converges, we can deduce the even/odd parts of  a monomial expansion in 1D using the even/odd monomials:
\[
\underbrace{\sum_{k=0}^\infty f_k x^k}_{\addition{f(x)}} \qquad = \qquad\underbrace{\sum_{k=0}^\infty f_{2k} x^{2k}}_{\addition{\fe(x)}}\qquad + \qquad  \underbrace{\sum_{k=0}^\infty f_{2k+1} x^{2k+1}}_{\addition{\fo(x)}}.
\]
We can write this in terms of the cyclic group of order 2, $\rmC_2 := \{-1,1\}$, capturing the reflection group action.  For any $g \in \rmC_2$ the even and odd monomials turn the symmetry action into multiplication by $\pm 1$:
\[
(g x)^{2k} = \underbrace{1}_{\sigmat(g)}   x^{2k}, \qquad (g x)^{2k+1} = \underbrace{g}_{\sigmas(g)} x^{2k+1}.
\]
Here $\sigmat,\sigmas : \rmC_2 \rightarrow O(1)$  are examples of orthogonal representations of $\rmC_2$, where  $O(n):= \{Q \in \bbR^{n \times n} : Q^\top Q = I\}$ denotes the set of orthogonal matrices and, in a slight abuse of notation, we  use the convention that $O(1) := \{\pm 1\}$ is scalar. They are maps from the group to orthogonal matrices (or, in this case, scalars) that preserve the group structure: for $a,b \in \rmC_2$ we have $\sigmat(ab) = \sigmat(a) \sigmat(b)$ and $\sigmas(ab) = \sigmas(a) \sigmas(b)$.   Using a monomial basis (or any basis that captures even/odd symmetry) in the discretisation of an ODE which commutes with reflection will automatically split into two linear systems according to parity.
}

\addition{
Two-dimensional monomials $x^k y^j$ similarly give a symmetry-adapted basis for horizontal/vertical reflection symmetry which corresponds to the product group $\rmC_2^2 := \{ \diag(a, b) : a,b \in \pm 1 \} \subset O(2)$.   In particular, for $g = \diag(a,b) \in \rmC_2^2$ reflections turns into multiplication by $1$ or $-1$ in four possible ways:
\meeq{
(a x)^{2k} (b y)^{2j} = \underbrace{1}_{\sigmat(g)} x^{2k} y^{2j},  &  (a x)^{2k+1} (b y)^{2j} &= \underbrace{a}_{\sigmah(g)} x^{2k+1} y^{2j}, \ccr
(a x)^{2k} (b y)^{2j+1} = \underbrace{b}_{\sigmav(g)} x^{2k} y^{2j+1},  &  (a x)^{2k+1} (b y)^{2j+1} &= \underbrace{a b}_{\sigmas(g)} x^{2k+1} y^{2j+1}.
}
Here $\sigmat, \sigmav, \sigmah, \sigmas : \rmC_2^2 \rightarrow O(1)$ are four examples of scalar orthogonal representations of $\rmC_2^2$ which all preserve the group structure.
Using a monomial basis (or any basis that captures even/odd symmetry) in the discretisation of a PDE that commutes with horizontal/vertical reflection will automatically decompose the PDE according into 4 independent linear systems. 
}

\addition{
How can we build a basis that captures all the symmetries of the square, i.e., the $D_4$ group? The set of homogeneous polynomials of degree $n$  are closed under the $D_4$ symmetry action so we first consider the monomial basis, which we put into a column vector:
\[
\bfp_n(x,y) := \vectt[x^n, x^{n-1} y, \ldots, x y^{n-1}, y^n].
\]
The symmetry actions of vertical reflection and $90^\circ$ rotation (which generate the group $D_4$) are now turned into matrix multiplication:
\meeq{
\bfp_n(x,-y) =  \sopmatrix{1 \\ & -1 \\ && \ddots \\ &&&  (-1)^n} \bfp_n(x,y), \ccr
\bfp_n(-y,x) =  \sopmatrix{&&&  (-1)^n\\
&& \iddots \\
& -1 \\
1}\bfp_n(x,y)
}
These two matrices actually generate an example of a matrix representation (they satisfy the same rules as the generators of $D_4$), but this representation is not {\it irreducible}: it is not the simplest way that  reflections and  rotations  can be represented as matrices. Instead, we will look for a basis so that the induced matrices are as simple as possible: we want them to correspond to irreducible representations of $D_4$.  We shall see that there are 5 irreducible representations; 4 are scalar whilst the faithful representation is two-dimensional. Building a symmetry-adapted basis  will decouple a PDE that commutes with reflection and  rotations across 6 independent linear systems, the number of which is equal to the sum of the dimensions of all the irreducible representations.
}

\addition{
This procedure can be generalised to the symmetry group of the cube, the {\it octohedral group} $O_h$. This group has 10 irreducible representations, mostly non-scalar, with a sum of the dimensions equal to 20. Rather than an explicit symmetry-adapted basis we construct such a basis numerically, thereby decoupling PDEs across 20 independent systems.  The techniques work for other groups, in particular the symmetric group $S_n$ corresponding to permutations with a more substantial  reduction in the dimension of the worst-case system when compared with monomials. Both the number of independent systems and the reduction in the dimension of the worst-case system increase rapidly with the dimension of PDEs with permutation symmetry. 
}

In the discussion we focus on monomial bases for pedagogical reasons but essentially everything we discuss is applicable to multivariate orthogonal polynomials, which have a natural link to symmetry groups in that the \addition{span of} degree-$n$ orthogonal polynomials are closed under \addition{orthogonal} symmetry actions \cite[Theorem 3.2.14]{DunklXu} and hence generate representations.  It should be emphasised for the reader inexperienced in numerical methods: \addition{{\it  avoid using monomials} for numerics}! They are prone to ill-conditioning and orthogonal polynomials have many beautiful properties such as sparse discretisations (as discussed in detail in  our recent review  \cite{SOActa} and references therein).  In fact many of the figures in this paper are generated  using multivariate orthogonal polynomials  for numerical stability \addition{reasons}.

The paper is structured as follows:

\noindent {\it \secref{background}}: We discuss prior work on incorporating discrete symmetries into the discretisation of partial differential equations via symmetry-adapted bases. \addition{This has primarily focussed on two-dimensional PDEs discretised using the Finite Element Method on regular 2D polygons (dihedral symmetry), whereas the present paper focusses on global spectral methods. We contend that the potential for these ideas having an impact are substantially higher in the spectral method setting, particularly as the dimension increases.}

\noindent {\it \secref{crash}}: We give a quick overview of the essentials of representation theory with an emphasis on practical applied mathematics: we state the standard results in terms of concrete matrices and vectors as opposed to abstract vector spaces. This includes {\it Schur's lemma} (\lmref{SchurMatrix}) which is a fundamental result in representation theory that states \addition{that a matrix which {\it intertwines} with two irreducible representations (multiplication on the left by one is the same as multiplication on the right by the other) is by necessity very sparse. This  will guarantee that discretisations of PDEs decouple for symmetry-adapted bases}.  It also includes a discussion on how any representation can be decomposed into irreducible representations via an orthogonal matrix which can be computed numerically.  

\noindent {\it \secref{symdecompos}}: Symmetry-adapted bases are constructed on a square  in closed form and on a cube \addition{numerically, utilising a} numerical algorithm for decomposing representations. \addition{An implication of this construction is a systematic approach for decomposing functions according to their symmetry classes, a la \figref{squaredecomposition}. }

\noindent {\it \secref{PDEs}}: Symmetry-adapted bases  decouple PDEs into distinct systems, whose number is the sum of the dimensions of the irreducible representations. The workhorse for proving this result is Schur's lemma: by showing that discretisations of PDEs are intertwining matrices we guarantee that they are sparse.  We demonstrate examples on the square using the dihedral group $D_4$ (which decouples across 6 systems) and a cube using the octohedral group $O_h$ (which decouples across 20 systems). 

\addition{\noindent {\it \secref{PermutationPDE}}: } \addition{We consider  using these techniques for   Schr\"odinger equations with multiple particles, beginning}
  with an example of 3 one-dimensional particles (so also living in a cube) using only permutation and negation symmetry (which decouples across 8 systems, a four-fold increase over \addition{standard} bases). \addition{For higher-dimensional Schr\"odinger equations we can deduce the number of independent systems and the reduction in the dimensions of the worst-case system, demonstrating a notable theoretical increase in computational effectiveness. Of particular interest are solutions with special symmetry properties corresponding to Bosons and Fermions where the reduction in dimensions is (super-)exponential.  }

\noindent {\it \secref{conc}}: We conclude by discussing extensions to vector-valued PDEs \addition{like Maxwell's equation} and  \addition{the potential for incorporation into preconditioners for PDEs that do not have the required symmetry properties to directly take advantage of representation theory in their discretisation.}

\Section{background} Background.

\addition{Whilst the idea of using symmetries in geometric integration gained prominence in the late 1980s \cite{kang1989construction,lasagni1988canonical,sanz1988runge}, the lesser known concept of symmetry-adapted bases incorporating knowledge from representation theory actually predates this field, with the first usage of dihedral symmetry for the numerical solution of PDEs on regular polygons via  Finite Element Methods appearing in the mid 1970s in the thesis of F\"assler \cite{fassler1976application} (see also \cite{stiefel2012group}\addition{, and the recent Python package implementing this approach for building symmetry-adapted bases \cite{PySymmetry}}).   Further investigation of tackling subproblems associated with such a construction was investigated by Bossavit  \cite{bossavit1993boundary}. Allgower, Georg and collaborators pursued using symmetry-adapted bases for Boundary Element Methods \cite{aimi2008restriction,allgower1998bem,allgower1998integral,allgower1998numerical}  including on the surface of a cube  \cite{georg1993exploiting}.
}
 
 \addition{While the existing work on incorporating representation theory into mesh-based methods (Finite and Boundary Element Methods) was mathematically beautiful there was no noticeable uptake in practical computation.  There are a couple reasons for this. There are very few geometries in 3D that can be approximated by meshes that have symmetries, essentially limited to  variants of the five platonic solids, or put another way,  any finite subgroup of 3D rotations ($SO(3)$) has only one of five forms \cite{aahlander2005applications}. But the real deal-killer is that    multigrid techniques enable significantly larger scale parallelisation than symmetry-adapted bases. Thus research on this topic largely fizzled out in the early 2000s.}

 
\addition{But symmetry is an inherently global phenomena, and therefore the application to global spectral methods is a much more natural avenue for these ideas. Further, the lack of an analogue of multigrid for spectral methods implies that there is significantly more potential for the ideas to have a material impact on the parallelisability of spectral methods. The  underlying idea of using symmetry and representation theory to decouple PDEs in a sense underlies the success of Fourier-based spectral methods (which capture translation symmetry)  as well as techniques based on spherical harmonics (which capture rotational symmetry, cf.~\cite{vasil2019tensor} and references therein). But the utilisation of discrete symmetry groups in spectral methods is highly underdeveloped with a small number of notable exceptions.    }  The $D_4$ symmetry group has  been used in the construction of spectral methods \cite{li2014symmetrizing}, \addition{using a construction that is equivalent to what is developed in \secref{D4basis}. Unfortunately, as we shall see, the benefits of the techniques for partial differential equations on a square are marginal as the resulting systems are no smaller than the standard techniques.} Munthe-Kaas constructed Fourier-like symmetry-adapted  bases associated with different discrete groups including the $D_3$ symmetry group associated with an equilateral triangle \cite{munthe2006group,StructurePreserving}.  \addition{But there does not appear to be any investigation of spectral methods on cubes or utilising permutation symmetry, which we shall see are the cases where the techniques may have the most profound impact, particularly in the high-dimensional setting.}

\addition{ Ahlander and Munthe-Kaas advocate in \cite{aahlander2005applications} an attractive  approach of incorporating symmetries directly on the level of linear algebra by using a group discrete Fourier transform to block-diagonalise the resulting linear systems, an approach that applies equally well for Finite Element and spectral method discretisations.  They present this idea via a beautiful analogy to using the discrete Fourier transform to diagonalise discretisations which capture translation symmetry. However, the number of independent systems in this approach is only equal to the number of irreducible representations, which may be significantly less than  the sum of their dimensions. Moreover, the computational complexity of the generalised Fourier transform grows with the order of the group (the number of elements)  which increases combinatorially fast for high-dimensional symmetry groups, limiting its applicability.
}


\addition{Beyond discretisation of PDEs, the underlying idea of using symmetry and representation theory to decouple PDEs in a sense underlies the success of  much of quantum mechanics, as discussed in  \cite{tinkham2003group,cotton1991chemical}.}
Outside of PDEs, symmetry-adapted bases have been used effectively in computational algebraic geometry  \cite{gatermann2004symmetry}, computing cubature rules \cite{HubertCubature}, and optimisation \cite{metzlaff2025symmetry}. \addition{An exciting recent development is in machine learning,via Equivariant Neural Networks \cite{cohen2016group} and Steerable Convolutional Neural Networks \cite{cohen2016steerable}, which have been used recently in the solution of inverse problems \cite{celledoni2021equivariant}. These  use notions of equivariance and symmetry-adapted bases to construct neural networks which capture symmetries.}


To extend these tools  \addition{beyond dihedral groups} requires the algorithmic decomposition of representations.  Dixon \cite{DixonComputeReps} introduced an iterative approach to decompose a representation  but the algorithm does not split apart multiple copies of the same irreducible representation.  Serre \cite{serre1977linear} gives an explicit construction of projectors including \addition{resolving} the case of multiple copies of the same irreducible representation.
This involved summing over every element of the group but the cost can be massively reduced   \cite{hymabaccus2020decomposing}, a technique which has been implemented in a GAP package  \cite{hymabaccus2020repndecomp}. An important special case that arises below are monomial representations, which have exactly one nonzero entry in each row and column, and an algorithm specific to \addition{this setting} was constructed in \cite{puschel2002decomposing}.  The method for incorporating symmetries of the cube in \cite{georg1993exploiting} utilised software developed in \cite{walker28numerical}, which is unfortunately no longer available. More recently, I introduced an algorithm based purely on numerical linear algebra (and implementable with standard floating point arithmetic)  for block-diagonalising representations of the symmetric group  \cite{SOSymGroup}.  We will utilise its numerical implementation  \cite{NumericalRepresentationTheory} to compute symmetry-adapted bases for the cube.

\Section{crash} A crash course in representation theory.

A representation is a map from a group to invertible linear operators in a way that the group structure is preserved. \addition{We focus on the specific case where the representations are orthogonal}:

\Definition{representation} \addition{For a group $G$ an {\it orthogonal representation} is a map $\rho : G \rightarrow O(n)$ that is a homomorphism, that is, for any $a,b \in G$ we have
\[
\rho(ab) = \rho(a) \rho(b).
\]}

\begin{table}[tb]\label{Table:Permutations}
\centering
\caption{\addition{Three equivalent ways of describing permutations of \( S_3 \). In cycle notation, a permutation maps each element in a cycle to the next element, and the last element of a cycle is mapped to the first element. In Cauchy's notation,  a permutation maps the top row to the bottom row. A permutation matrix acts on a vector by standard matrix-vector multiplication,  permuting its rows.  For concreteness, we treat all permutations in $S_3$ as equivalent to their permutation matrices.} }
\[
\begin{array}{c c c}
\hline 
\hbox{Cycle Notation} & \hbox{Cauchy's Notation} & \hbox{Permutation Matrix} \\ 
\hline   \\[-2mm]
 (1)(2)(3)  &  \begin{pmatrix} 1 & 2 & 3 \\ 1 & 2 & 3 \end{pmatrix}  &  \begin{pmatrix} 1 & 0 & 0 \\ 0 & 1 & 0 \\ 0 & 0 & 1 \end{pmatrix}  \\[6mm]
 (1\ 2)(3)  &  \begin{pmatrix} 1 & 2 & 3 \\ 2 & 1 & 3 \end{pmatrix}  &  \begin{pmatrix} 0 & 1 & 0 \\ 1 & 0 & 0 \\ 0 & 0 & 1 \end{pmatrix}  \\[6mm]
 (1\ 3)(2)  &  \begin{pmatrix} 1 & 2 & 3 \\ 3 & 2 & 1 \end{pmatrix}  &  \begin{pmatrix} 0 & 0 & 1 \\ 0 & 1 & 0 \\ 1 & 0 & 0 \end{pmatrix}  \\[6mm]
 (1)(2\ 3)  &  \begin{pmatrix} 1 & 2 & 3 \\ 1 & 3 & 2 \end{pmatrix}  &  \begin{pmatrix} 1 & 0 & 0 \\ 0 & 0 & 1 \\ 0 & 1 & 0 \end{pmatrix}  \\[6mm]
 (1\ 2\ 3)  &  \begin{pmatrix} 1 & 2 & 3 \\ 2 & 3 & 1 \end{pmatrix}  &  \begin{pmatrix} 0 & 0 & 1 \\ 1 & 0 & 0 \\ 0 & 1 & 0 \end{pmatrix}   \\[6mm]
 (1\ 3\ 2)  &  \begin{pmatrix} 1 & 2 & 3 \\ 3 & 1 & 2 \end{pmatrix}  &  \begin{pmatrix} 0 & 1 & 0 \\ 0 & 0 & 1 \\ 1 & 0 & 0 \end{pmatrix}  \\[6mm]
\hline
\end{array}
\]
\end{table}

As a  basic example consider the symmetric group $S_3$ corresponding to permutations of 3 elements. \addition{ There are a few standard  notations for  permutations---Cauchy's notation, cycle notation,  and permutation matrices\footnote{\addition{Our convention for permutation matrices is the transpose of the standard convention in order to preserve the group structure.}}---which we depict in \tabref{Permutations} for $S_3$. For concreteness we define the group $S_3$ to be the set of permutation matrices, which can be viewed as a symmetry group of the unit cube: permuting variables leaves the cube unchanged. The permutation group $S_3$ is generated from the simple transpositions
\[
\tau_1 := (1\ 2) (3) = \begin{pmatrix} 0 & 1 \\ 1 & 0 \\ && 1 \end{pmatrix}\qqand \tau_2 := (1) (2\ 3) = \begin{pmatrix} 1\\ &0 & 1 \\ &1 & 0 \end{pmatrix},
\] 
that is, every element of $S_3$ can be written as a product of $\tau_1$ and $\tau_2$, which we can write as  $S_3 = \ip<\tau_1,\tau_2>$.  For example, the permutation $(1\ 2\ 3)$ is a product of two simple permutations since:
\[
\tau_1 \tau_2 = \sopmatrix{0 & 0 &1  \\  1 & 0 & 0 \\ 0 & 1 & 0} =  (1\ 2\ 3).
\]
}

\begin{table}[tb]
\centering
\caption{\addition{ Canonical (orthogonal) irreducible representations of \( S_3 \). The traces of this table give the classic character table for $S_3$, but we will not use characters. The whole table can in fact be deduced from any generating set such as the simple transpositions, i.e., the second and fourth row. }}\label{Table:Irreducible}
\begin{tabular}{c | c c c}
\hline\\[-3mm]
{Cycle Notation} & {$\sigma_3$ (Trivial)} & {$\sigma_{2+1}$} & {$\sigma_{1+1+1}$ (Sign)} \\[1mm]
\hline\\[-2.5mm]
\( (1)(2)(3) \) & \( 1 \) & \( \begin{pmatrix} 1 & 0 \\ 0 & 1 \end{pmatrix} \) & \( 1 \) \\[4mm]
\( (1\ 2)(3) \) & \( 1 \) & \( \begin{pmatrix} -1 &0  \\ 0 & 1 \end{pmatrix} \) & \( -1 \) \\[4mm]
\( (1\ 3)(2) \) & \( 1 \) & \( {1 \over 2} \begin{pmatrix} 1 & -\sqrt{3} \\ -\sqrt{3} &-1 \end{pmatrix} \) & \( -1 \) \\[4mm]
\( (1)(2\ 3) \) & \( 1 \) & \( {1 \over 2} \begin{pmatrix} 1 & \sqrt{3} \\ \sqrt{3} &-1 \end{pmatrix} \)  & \( -1 \) \\[4mm]
\( (1\ 2\ 3) \) & \( 1 \) &  \( {1 \over 2} \begin{pmatrix} -1 & -\sqrt{3} \\ \sqrt{3} &-1 \end{pmatrix} \) & \( 1 \) \\[4mm]
\( (1\ 3\ 2) \) & \( 1 \) &  \( {1 \over 2} \begin{pmatrix} -1 & \sqrt{3} \\ -\sqrt{3} &-1 \end{pmatrix} \)  & \( 1 \) \\[4mm]
\hline
\end{tabular}
\end{table}

\addition{Permutation matrices themselves are the most natural example of an orthogonal representation of $S_3$:  the {\it standard representation} $\rhos : S_3 \rightarrow O(3)$ defined by $\rhos(g) = g$  is trivially a group homomorphism. But  this is not the only way of representing the permutation group as a set of matrices and in \tabref{Irreducible} we give three basic examples: two scalar representations $\sigma_3, \sigma_{1+1+1} : S_3 \rightarrow O(1)$ and one $2 \times 2$ representation $\sigma_{2+1} : S_3 \rightarrow O(2)$.  We will see in the next section that these are the {\it irreducible representations}, which will also explain the notation. We invite the reader to confirm these are indeed group homomorphisms.  For example, we can see how  the formula $\tau_1 \tau_2 =   (1\ 2\ 3)$ is preserved by each of these representations:
\meeq{
\underbrace{ 1 } _{\sigma_{3}\pr( \tau_1 \tau_2 )} =  \underbrace{ 1  } _{\sigma_{3}\pr( \tau_1 )}  \times  \underbrace{  1} _{\sigma_{3}\pr( \tau_2 )}, \ccr
\underbrace{  {1 \over 2} \begin{pmatrix} -1 & -\sqrt{3} \\ \sqrt{3} &-1 \end{pmatrix}  } _{\sigma_{2+1}\pr(\tau_1 \tau_2)} =  \underbrace{ \begin{pmatrix} -1 & 0\\ 0 & 1 \end{pmatrix}  } _{\sigma_{2+1}\pr( \tau_1 )}  \times  \underbrace{ {1 \over 2} \begin{pmatrix} 1 & \sqrt{3} \\ \sqrt{3} &-1 \end{pmatrix} } _{\sigma_{2+1}\pr( \tau_2)}, \ccr
\underbrace{ 1 } _{\sigma_{3}\pr( \tau_1 \tau_2 )} =  \underbrace{ (-1)  } _{\sigma_{3}\pr( \tau_1 )}  \times  \underbrace{  (-1) } _{\sigma_{3}\pr( \tau_2 )}.
}
We have thus seen that $S_3$ has at least four different representations: the standard permutation matrices and three irreducible representations.  But we shall see that permutation matrices can in fact be reduced to two of these irreducible representations. 
}

\addition{
\begin{remark}
The reader familiar with group theory may have seen the concept of a {\it character table}, which is defined as the trace of the irreducible representations. Character tables play a fundamental role in group theory, for example, they are used in the definition of the group discrete Fourier transform \cite{aahlander2005applications}, but will not be necessary below.  Whereas a character table requires knowledge of all elements of the group, the table of irreducible representations can in fact be deduced from only knowledge of any generating set: if  $g_1,\ldots,g_m$ are  generators of $G$ then any element of the group can be written as a product  ($g = g_{k_1} \cdots g_{k_M}$)
and the representation deduced by multiplying the representations of the generators accordingly ($
\rho(g) = \rho(g_{k_1}) \cdots \rho(g_{k_M})$).
Thus in \tabref{Irreducible} we in fact only need to know the second and fourth row in order to deduce the entire table.
\end{remark}
}


\Subsection Reducible and irreducible representations.



	
\addition{ A subspace $V \subset \bbR^n$  is invariant under a representation $\rho$ if it is closed under multiplication by the image of $\rho$: for all $\bfv \in V$ and $g \in G$ we have $\rho(g) \bfv \in V$. An orthogonal representation is {\it reducible} if there exists a {\it nontrivial}  invariant subspace $V$, i.e., $V$ is not empty or the full space $\bbR^n$.   }For example,  \addition{the standard permutation matrices}  leave the space of vectors with all entries equal ($\{\vectt[c,\dots,c] : c \in \bbR \}$) invariant \addition{as permuting the entries does not modify the vector at all}, thus the standard representation $\rhos$ is in fact a {\it reducible representation}. \addition{An {\it irreducible representation} is a representation that is not reducible: the only invariant subspaces are the empty space and $\bbR^n$.}

Different irreducible representations can be {\it isomorphic}: \addition{if $\rho : G \rightarrow O(n)$ is an irreducible representation and $Q \in O(n)$ then $\tilde\rho(g) :=  Q^\top \rho(g) Q$ is also an irreducible representation, as if one has a non-trivial invariant subspace than so does the other}. However there are only a finite number of irreducible representations up to isomorphism:

\Theorem{numberirreps}\cite[Theorem 7]{serre1977linear}
	The number of irreducible representations of a group $G$ (up to isomorphism) is equal to the number of conjugacy classes of $G$, \addition{i.e., the number of equivalence classes under the relationship $a \sim b$ \addition{for $a,b \in G$} if there exists $g \in G$ such that $a = g b g^{-1}$}.

In other words given any finite group we can enumerate the irreducible representations, choosing a specific form for each isomorphism class. We refer to irreducible representations in such an enumeration as {\it \addition{canonical} irreducible representations}. \addition{We will typically denote canonical irreducible representations using $\sigma$ and all other representations using $\rho$.}

As a very simple example consider the cyclic group of order 2, $\addition{\rmC_2 := \{\pm 1\} = O(1)}$. Being a commutative group the \addition{number of} conjugacy classes are equal to the order \addition{of the group:  for all $g \in G$ we have $gag^{-1} = g g^{-1} a = a$, hence, any elements which are conjugate to each other must be the same. We thus} know there are \addition{precisely} two irreducible representations \addition{of $\rmC_2$}. These are in fact the trivial ($\addition{\sigmat(g) = 1}$) and sign representations ($\addition{\sigmas(g) = g}$), \addition{which  correspond to even and odd functions, respectively, when incorporated into a symmetry-adapted basis}.  A more complicated example is $S_3$  \addition{where two permutations are conjugate to each other if their cycles have the same lengths. Thus there are three conjugacy classes, which we identify with the three partitions of $3$ according to the length of their cycles when ordered from largest to smallest}:
\[
\underbrace{\{(123), (132)\}}_{3 = 3},\quad \underbrace{\{(12)(3), (23)(1), (13)(2)\}}_{3 = 2+1},\quad \underbrace{\{(1)(2)(3)\}}_{3 = 1+1+1}.
\]
\addition{A consequence of there being exactly three conjugacy classes is that there are exactly three irreducible representations (up to isomorphism), precisely those defined in \tabref{Irreducible}. These irreducible representations will tell us the ways that permutation symmetry in 3D can be built into a symmetry-adapted basis.}

\begin{remark}
Note whilst the  irreducible representations and the conjugacy classes are both one-to-one with partitions it is a mistake to view them as some how connected to each other. The proof of \thref{numberirreps} is based on showing that they \addition{both lead to} bases of class functions (functions that depend only on conjugacy class) and hence have the same number. 
\end{remark}

\addition{
\begin{remark}
A fundamental fact is that the sum of the  {\it squares} of the dimensions of the irreducible representations is equal to the order of the group, e.g. in the case of the symmetric group $S_n$ it equals $n!$. However, when we decouple PDEs the number of independent systems equals  the sum of the dimensions, not their square, which does not have a simple formula in general.
\end{remark}
}

We now go through the enumeration of the irreducible representations for the square ($D_4$) and cube ($O_h$).

\Subsection Irreducible representations for the dihedral group $D_4$.

The symmetry group of an $n$-sided regular polygon is the {\it dihedral group}
$D_n$, and the general classification of irreducible representations is known  \cite[Section 5.3]{serre1977linear}.  Whilst the general construction has a simple  explicit form, for concreteness we specialise the results  to the  square $D_4$
which is generated by a rotation $r$ and reflection $s$:
\[
D_4 = \ip<r,s : r^4 = s^2 = (s r)^2 = 1> = \{I,\ \underbrace{r,r^2,r^3}_{\hbox{rotations}},\ \underbrace{s, sr, sr^2, sr^3}_{\hbox{reflections}} \}.
\]
\addition{To minimise notation we will identity $D_4 \subset O(2)$, where vertical reflection and a $90^{\circ}$ counter-clockwise rotation take the explicit form:
\[
s = \begin{pmatrix} 1  & 0\\ 0& -1 \end{pmatrix}\qand r = \begin{pmatrix} 0 & -1 \\ 1 & 0 \end{pmatrix}.
\]
Thus $D_4$ can be viewed as a symmetry group which leaves  the unit square  $[-1,1]^2$ invariant.  }

\addition{One can verify that there are precisely
  five conjugacy classes:}
\[
\set{I}, \set{r,r^3}, \set{r^2}, \set{s ,  s r^2}, \set{sr, sr^3}. 
\]
\addition{and therefore we aim} to construct a total of 5 irreducible representations, \addition{which we define in \tabref{DihedralIrreducible}.  We will see in \secref{symdecompos} how these irreducible representations can be built into symmetry-adapted bases in order to compute the symmetry decomposition in \figref{squaredecomposition}, where the faithful representation will tell us precisely how the symmetry is present in $f_{\rm f1}$ and $f_{\rm f2}$.}
 



\begin{table}[tb]
\centering
\caption{\addition{Irreducible representations of the Dihedral Group $D_4$. The entire table can be deduced from the generators, given in the second and fifth row. We will construct symmetry-adapted bases where reflections and rotations are equivalent to multiplications by one of these irreducible representations.}} \label{Table:DihedralIrreducible}
\[\begin{array}{c|ccccc}
\hline\\[-2.5mm]
\hbox{Element} & \sigmat\ \hbox{(Trivial)} & \sigmats\ \hbox{(Reflect)} & \sigmast\ \hbox{(Rotate)} & \sigmas\ \hbox{(Sign)} & \sigmaf\ \hbox{(Faithful)} \\[1mm]
\hline\\[-2.5mm]
I     & 1  &  1  &  1  &  1  &  \sopmatrix{1 & 0 \\ 0 & 1}  \\[3mm]
r     & 1  &  1  & -1 & -1 &  \sopmatrix{0 & -1 \\ 1 & 0}    \\[3mm]
r^2   & 1  &  1  &  1  &  1  &  \sopmatrix{ -1 & 0 \\0 & -1}  \\[3mm]
r^3   & 1  &  1  & -1 & -1 &  \sopmatrix{ 0 & 1  \\-1 & 0}   \\[3mm]
s   & 1  & -1 &  1  & -1 &   \sopmatrix{1 & 0 \\ 0 & -1}   \\[3mm]
sr   & 1  & -1 &  -1  & 1 &   \sopmatrix{0 & -1 \\ -1 & 0}   \\[3mm]
sr^2   & 1  & -1 & 1 &  -1  &   \sopmatrix{ -1 & 0 \\0 & 1}    \\[3mm]
sr^3   & 1  & -1 & -1 &  1  &  \sopmatrix{ 0 & 1  \\1 & 0}   \\[3mm]
\hline
\end{array}\]
\end{table}

\subsection{Irreducible representations for the octohedral group $O_h$}\label{Section:cubeirreducibles}

\addition{Hyperoctohedral groups are  symmetry groups for hypercubes for which the irreducible representations are known in general \cite{geissinger1978representations,musili1993representations}}. In the special case of the octohedral group \addition{for the 3D cube} they can be connected directly to that of the symmetric group $S_4$, \addition{that is, the group of permutations of 4 elements}. In particular, we have $O_h \cong S_4 \times \addition{\rmC}_2$: the cube has four diagonals and any symmetry action permutes these four diagonals, alongside a choice of reflection, see \figref{cubediagonals}.  We  have three generators \addition{$\tau_1,\tau_2,\tau_3$} corresponding to permuting the diagonals and one generator $s$  corresponding to negation. 
Explicitly,  we can write the four diagonals as
\[
\bfd_1 =  \begin{pmatrix} 1\\1\\1 \end{pmatrix} \bbR,\quad \bfd_2 =  \begin{pmatrix} 1\\-1\\1 \end{pmatrix} \bbR,\quad \bfd_3 =  \begin{pmatrix} 1\\1\\-1 \end{pmatrix} \bbR,\quad \bfd_4 =   \begin{pmatrix} 1\\-1\\-1 \end{pmatrix} \bbR.
\]
We see that the simple transpositions are given by
\[
\tau_1 = \begin{pmatrix}0 &0& 1 \\ 0& -1 & 0 \\ 1 & 0 & 0 \end{pmatrix}, 
\tau_2 = \begin{pmatrix}  -1 & 0 & 0 \\0 &0& -1 \\ 0 & \addition{-}1 & 0 \end{pmatrix}, \tau_3 = \begin{pmatrix} 0&0& -1 \\ 0& -1 & 0 \\ -1 & 0 & 0 \end{pmatrix},
\]
\addition{which} permute the diagonals: for $k=1,2,3$ we have
$
\tau_k \bfd_k = \bfd_{k+1}, \tau_k \bfd_{k+1} = \bfd_k
$
and all other diagonals remain fixed. The last generator of $O_h$ is negation:
$
s = -I. 
$

\Figurefixed[cubediagonals]
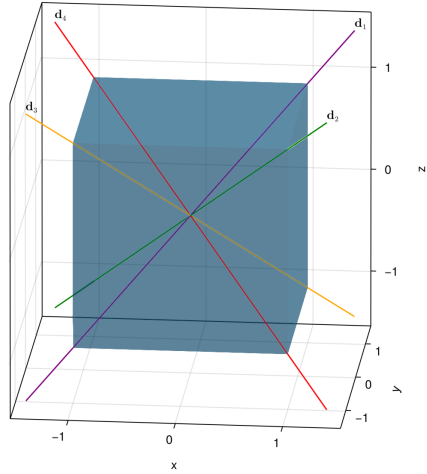
	The four diagonals of a cube. Rotational symmetries of the cube correspond to permuting these diagonals. These 24 permutations combined \addition{with} negation define the 48 symmetries of the cube $O_h \cong S_4 \times \addition{\rmC}_2$. 

A recipe exists for constructing the \addition{simple transposition} generators \addition{$\tau_1,\ldots,\tau_{n-1}$} of all orthogonal irreducible representations of the symmetric group \cite{Okounkov}. In particular, as in the case of $S_3$ the conjugacy classes of $S_n$ correspond to cycles with specified lengths and hence irreducible representations can be put in one-to-one correspondence with  partitions of $n$.   In the case of $S_4$, for the \addition{five} different partitions $\lambda$  \addition{of 4}  the recipe yields the  formulae \addition{in \tabref{OctohedralIrreducible}} for the generators $\sigma_\lambda(\tau_k)$, \addition{which suffices to define the representation by representing any element of $O_h$ as a product of the generators (including the whole table with 24 columns would be excessive!)}.


\begin{table}[tb]
\centering
\caption{\addition{Generators for the irreducible representations of the symmetric group $S_4$, where $\lambda$ is a partition of 4.  These irreducible representations are used to construct the irreducible representations of the octohedral group $O_h$ corresponding to the symmetries of the cube.} \label{Table:OctohedralIrreducible}}
\[
\begin{array}{ c|c|c| c } 
\hbox{Partition $\lambda$} & \sigma_\lambda(\tau_1) & \sigma_\lambda(\tau_2) & \sigma_\lambda(\tau_3) \\[5pt]
 \hline
&&&  \\[-8pt]
4 & 1 & 1 & 1 \\[3pt] 
 3+1 & \smallpmatrix{  -1 \\ & 1 \\ & & 1  }&   \smallpmatrix{ 1/2 & \sqrt{3}/2  \\  \sqrt{3}/2  & -1/2  \\ && 1} & \smallpmatrix{  1 \\ &1/3 & 2\sqrt{2}/3 \\  &2\sqrt{2}/3 &  -1/3 &} \\[12pt] 
2+2 &   \smallpmatrix{  -1 \\ & 1   }  &  \smallpmatrix{ 1/2 & \sqrt{3}/2  \\  \sqrt{3}/2  & -1/2 } &   \smallpmatrix{  -1 \\ & 1   } \\[12pt] 
 2+1+1	 & \smallpmatrix{  -1 \\ & -1 \\ & & 1  }  & \smallpmatrix{  -1 \\ &  1/2 & \sqrt{3}/2 \\ &   \sqrt{3}/2 & -1/2  } &  \smallpmatrix{  1/3 & 2\sqrt{2}/3 \\  2\sqrt{2}/3 &  -1/3 &  \\ &    & -1  } \\[12pt]
1+1+1+1 & -1 & -1 & -1
\end{array}
\]
\end{table}

The irreducible representations of a product group are tensor products of the irreducible representations of each component \cite[Section 3.2]{serre1977linear}, which since the irreducible representations of $\addition{\rmC}_2$ are scalar are equivalent to a standard product. That is, we have $10$ irreducible representations which we define for any of the five partitions $\lambda$ via their generators as
\[
\sigma_{\lambda,{\rm t}/{\rm s}}(\tau_k) := \sigma_\lambda(\tau_k), \qquad \sigma_{\lambda,{\rm t}}(s) = I, \qquad  \sigma_{\lambda,{\rm s}}(s) = -I.
\]

\addition{How does the standard representation of the octohedral group generated by $\tau_1, \tau_2, \tau_3$ and $s$ correspond to the irreducible representations of $O_h$? The representation  $\rho(g) = g$},   is irreducible but is not canonical: it does not appear in \addition{\tabref{OctohedralIrreducible}}. However, we can reduce it explicitly to canonical form using the matrix
\[
Q = \sopmatrix{1/\sqrt{3} &1/\sqrt{6} & -1/\sqrt{2} \\ -1/\sqrt{3} & 2/\sqrt{6} & 0 \\ -1/\sqrt{3} & -1/\sqrt{6} & -1/\sqrt{2}}.
\]
So that $Q^\top \rho(g) Q = \sigma_{2+1+1,{\rm s}}(g)$, which can be verified by checking the generators.

\Subsection Schur's lemma.

Whilst representation theory is a rich and beautiful theory for our primary purpose of showing that a PDE splits across different symmetry classes we only need one foundational result: {\it Schur's lemma}.  This is a powerful tool that tells us something about so-called {\it intertwining maps}. \addition{For simplicity we rephrase this concept in terms of the equivalent concept of intertwining matrices:}

\Definition{intertwine}
\addition{ Let $\rho_1 : G \rightarrow O(m)$ and $\rho_2 : G \rightarrow  O(n)$ be two orthogonal representations. An intertwining matrix $A \in \bbR^{m \times n}$ with respect to $\rho_1$ and $\rho_2$  satisfies, for all  $g \in G$,
\[
\rho_1(g) A = A \rho_2(g).
\] 
}

Schur's lemma (see e.g. \cite[Schur's Lemma 1.7]{FultonRepTheory}) states that  intertwining maps are highly structured.  
Hidden inside \addition{the usually abstractly stated lemma} is an elegant tool for proving that specific matrices are either zero or a constant times the identity \addition{(i.e., {\it very sparse} matrices!)}: 

\begin{lemma}[\addition{Schur's lemma, rephrased.}]\label{Lemma:SchurMatrix}
\addition{Suppose $\sigma_1 : G \rightarrow O(m)$ and $\sigma_2  : G \rightarrow O(n)$ are  canonical irreducible representations. Suppose that  $A \in \bbR^{m \times n}$  is an intertwining matrix with respect to $\sigma_1$ and $\sigma_2$. For some $\lambda \in  \bbR$ we have
\[
A = \begin{cases} 0 & \sigma_1 \neq \sigma_2, \\
			\lambda I & \sigma_1 = \sigma_2.
			\end{cases}
\]}
\end{lemma}

\begin{proof}
\addition{As it is such a fundamental result to this paper and it is not quite in  standard form, we reproduce the proof of Schur's lemma here.  
The kernel $K$ of $A$ is invariant under multiplication by $\sigma_2$: if $\bfk \in K$ then for all $g \in G$ we have $\sigma_2(g) \bfk \in K$ since
\[
 A \sigma_2(g) \bfk = \sigma_1(g) A  \bfk = 0.
\]
By the definition of irreducibility we therefore know that $K$ is  either empty or is equal to $\bbR^n$, that is $A$ has column rank of either 0 or $n$.  Similarly, the co-kernel $C$ of $A$ is invariant under multiplication by $\sigma_1$ and thus  $A$ has row rank of either 0 or $m$. If $A$ has rank 0 then it is zero. Otherwise, $A$ must be square (since the row and column ranks are equal) and invertible. This implies that $\sigma_1$ and $\sigma_2$ are isomorphic:
$
\sigma_1(g) = A^{-1} \sigma_2(g) A.
$
Since we are assuming they are both in canonical form this implies that $\sigma_1 = \sigma_2$. 
}

\addition{
Let $\lambda \in \bbC$  be an eigenvalue of $A$ and let $V \subset \bbC^n$ denote the corresponding eigenspace. Note that $A- \lambda I $ is also an intertwining matrix whose kernel is $V$, and by the same logic as above $V$ must be invariant  under $\sigma_1 = \sigma_2$.  But since $\sigma_1$ is irreducible and $V$ is not empty we know that $V = \bbR^n$, that is every vector is in the kernel of $A-\lambda I$, and hence $A = \lambda I$.  Since $A$ has real entries we are guaranteed that $\lambda \in \bbR$.
}
\end{proof}

\addition{This sparsity property will prove invaluable for block-diagonalising discretisations of PDEs.}


%
%

\Subsection Decomposing  into irreducible representations.

An intuitive result is that a representation can be decomposed into subrepresentations until all we are left with are irreducible representations. 
\addition{D}ecomposing a representation is best thought of as a block-diagonalisation process akin to diagonalising a matrix.  We express this as follows:

\Theorem{fullblockdiag}\cite[Theorem 1]{SOSymGroup} Given an orthogonal representation $\rho : G \rightarrow O(n)$ of a finite group $G$, there exists an orthogonal matrix $Q$ such that
\[
Q^\top \rho(g) Q = \sopmatrix{ \addition{\sigma}_1(g)^{\oplus a_1} \\ & \ddots \\ && \addition{\sigma}_r(g)^{\oplus a_r}}
\]
where $\addition{\sigma}_1,\ldots,\addition{\sigma}_r$ are all \addition{canonical} irreducible representations.  Here we use the notation
\[
A^{\oplus k} := \diag(\underbrace{A,\ldots,A}_{\hbox{$k$ \rm times}}).
\]

As an example we consider block-diagonalising the permutation matrices associated with $S_3$. The matrix
\[
Q = \pr(\begin{array}{cc|c}
\sqrt{2} & 1/\sqrt{6} & -1/\sqrt{3} \\
-\sqrt{2} & 1/\sqrt{6} & -1/\sqrt{3} \\
 & -\sqrt{2/3} & -1/\sqrt{3}
\end{array})
\]
satisfies
\meeq{
Q^\top \begin{pmatrix} 0 & 1 \\ 1 & 0 \\ && 1 \end{pmatrix} Q = \pr({\begin{array}{cc|c} -1 && \\ & 1& \\\hline && 1\end{array}}), \ccr
 Q^\top \begin{pmatrix} 1 \\ & 0 & 1 \\ & 1 & 0  \end{pmatrix} Q =  \pr({\begin{array}{cc|c}  1/2 & \sqrt{3}/2 \\\sqrt{3}/2  & -1/2  \\\hline && 1\end{array}}).
}
In other words, we have block-diagonalised the representation associated with permutation matrices to two  irreducible representations of smaller dimension: 
\[
Q^\top \addition{\rhos(g)} Q = \diag(\sigma_{2+1}(g), \sigma_3(g)).
\]

A number of algorithms \addition{exist based on} group theoretical considerations (character tables, cosets, etc.) \cite{DixonComputeReps,serre1977linear,hymabaccus2020decomposing,puschel2002decomposing}. However, an algorithm   built using only  numerical linear algebra\footnote{\addition{This algorithm can be viewed as a generalisation of the approach used in \cite{cohen2016group} for building Equivariant Neural Networks, which computed equivariant bases by computing the columns of $Q$ associated with the trivial representation.}} was introduced in \cite{SOSymGroup} that can compute the matrix $Q$   directly from generators of the representation.  This was shown only for the symmetric group  though the result \cite[Theorem 2]{SOSymGroup} translates to other groups which we outline now. Essentially the idea mimics an approach for computing eigenvectors of a matrix $A \in \bbR^{n \times n}$ via a nullspace computation. That is,  if a (known) eigenvalue $\lambda_k$ has (unknown) multiplicity $a_k$, a basis of eigenvectors  $\bfu_k^{(1)}, \ldots, \bfu_k^{(a_k)} \in \bbR^n$ can be computed  by finding a basis for solutions to
\[
A \bfu_k^{(j)} = \lambda_k  \bfu_k^{(j)},
\]
\addition{which is equivalent to a nullspace calculation for $A - \lambda_k I$. This procedure generalises to representations:}

\begin{lemma}
\addition{If we know  irreducible representations  $\sigma_k : G \rightarrow O(d_k)$  of a group $G$ generated by $g_1,\ldots,g_m$ then an orthogonal matrix $Q \in O(n)$ that decomposes a representation $\rho : G \rightarrow O(n)$ as  in \thref{fullblockdiag} can be computed by finding a basis of $U_k^{(j)} \in \bbR^{n \times d_k}$ that satisfy
\meeq{
\rho(g_1)  U_k^{(j)} = U_k^{(j)} \sigma_k(g_1), \\
&\vdots \ccr
 \rho(g_m)  U_k^{(j)} = U_k^{(j)} \sigma_k(g_m).
}}
\end{lemma}

\begin{proof}
\addition{
Any orthogonal matrix that satisfies the coniditons of \thref{fullblockdiag} can be written as
 \[
 Q = \left( \begin{array} {c : c  : c | c | c : c  : c}
 Q_1^{(1)}& \ldots& Q_1^{(a_1)}  &\cdots &  Q_r^{(1)}& \ldots& Q_1^{(a_r)}
 \end{array} \right)
 \]
 where each $Q_k^{(j)}$ satisfies the same conditions as $U_k^{(j)}$ in the above system. We want to show the other way around:  if we compute $\{U_k^{(j)}\}$ we can recover $\{Q_k^{(j)}\}$ that form a suitable orthogonal matrix.  
Vectorising the columns of $U_k^{(j)}$ reduces \addition{its determination} to a standard  nullspace problem.   Then $U_k^{(j)}$ is an intertwining matrix with respect to $\rho$ and $\sigma_k$: writing any element of the group as a product of generators $g = g_{\ell_1} \cdots g_{\ell_M}$ we find
\[
U_k^{(j)} \sigma_k(g) = U_k^{(j)} \sigma_k(g_{\ell_1}) \cdots \sigma_j(g_{\ell_M}) = \rho(g_{\ell_1}) \cdots \rho(g_{\ell_M})   U_k^{(j)}  = 
\rho(g)  U_k^{(j)}.
\]
Therefore  $(U_k^{(j)})^\top U_k^{(j)}$ is an intertwining matrix with respect to $\sigma_k$:
\[
(U_k^{(j)})^\top U_k^{(j)} \sigma_k(g) =  (U_k^{(j)})^\top \rho(g) U_k^{(j)}  =  (\rho(g^{-1})U_k^{(j)})^\top  U_k^{(j)} 
=  \sigma_k(g)   (U_k^{(j)})^\top  U_k^{(j)},
\]
where we use the fact for any orthogonal representation $\rho$ we have $\rho(g)^\top = \rho(g)^{-1} = \rho(g^{-1})$.
Schur's Lemma (\lmref{SchurMatrix}) ensures it is a (nonzero) constant times the identity.  We can thus rescale $U_k^{(j)}$ to determine $Q_k^{(j)}$ so that the inner products of its columns satisfy $(Q_k^{(j)})^\top Q_k^{(j)} = I$. Hence we can build an orthogonal $Q$ that block-decomposes $\rho$. 
}
\end{proof}

Note in all our examples we either construct   $Q$ explicitly or the representation is that of the symmetric group, hence the existing algorithm of \cite{SOSymGroup} and its implementation in \cite{NumericalRepresentationTheory} suffice.

\Section{symdecompos}  Computing symmetric decompositions.

We now turn to the task of computing symmetric decompositions such as those associated with $D_4$ used to produce \figref{squaredecomposition}.  The starting point is to assume that we have a basis of functions $\set{p_1,\ldots,p_N}$ where a symmetry action is equivalent to multiplying the basis by a representation. 
We write a basis as a column-vector:
\[
\bfp(\bfx) := \vectt[p_1(\bfx), \dots, p_N(\bfx)].
\]
Then we can express the notion of a symmetry-respecting and symmetry-adapted basis as follows:

\Definition{symrespect} A  {\rm  symmetry-respecting basis} with respect to an orthogonal representation $\rho : G \rightarrow O(N)$ is a vector of functions $\bfp : \Omega \rightarrow \bbR^N$ satisfying
\[
\bfp(g \bfx) =   \rho(g) \bfp(\bfx).
\]
If $\rho$ is a canonical irreducible representation we refer to it as an {\rm irreducible symmetry-adapted basis} and if $\rho$ is block-diagonal where each block is \addition{a canonical irreducible representation} (a la \thref{fullblockdiag}) then we refer to it as a {\rm symmetry-adapted basis}.

Our main example of a symmetry-respecting basis are  monomials of precisely degree $n$ which span homogeneous polynomials. For example, in 2D we have the basis
\[
\bfp_n(x,y) = \vectt[x^n, x^{n-1}y,\dots,x y^{n-1}, y^n] \in \bbR^{n+1}.
\]
We know (via the binomial theorem) that homogeneous polynomials are closed under any orthogonal change-of-variables, that is,  for any $g \in G \subset O(2)$ there must exist $\rho(g) \in O(n+1)$ such that
$
\bfp(g \bfx) =\rho(g) \bfp(\bfx)
$. It follows that $\rho$ must itself be a representation as we have
\[
\rho(ab)\bfp(\bfx) = \bfp(ab \bfx) = \rho(a)\bfp(b \bfx) = \rho(a)\rho(b) \bfp(\bfx) 
\]
hence $\rho(ab) = \rho(a) \rho(b)$.


\addition{W}e can employ \thref{fullblockdiag} \addition{(as realised in the algorithm of  \cite{SOSymGroup})} to deduce a matrix $Q$ that block-diagonalises \addition{$\rho$ and} thereby deduce a symmetry-adapted basis:

\Lemma{symrepOPs} 
Suppose $\bfp$ is a symmetry-respecting basis with respect to $\rho : G \rightarrow O(n)$ and $Q \in O(n)$ block-diagonalises $\rho$ \addition{into irreducible representations}. Then
\[
\bfq(\bfx) :=  Q^\top \bfp(\bfx)
\]
is a symmetry-adapted basis.

\Proof

This falls out immediately from the definitions:
\meeq{
\bfq(g \bfx) =  Q^\top \bfp(g \bfx)   =Q^\top \rho(g) \bfp(\bfx)  = \underbrace{Q^\top \rho(g) Q}_{\hbox{block-diagonalised}} \bfq(\bfx).
}
\mqed


\subsection{Symmetry-adapted basis with respect to $D_4$} \label{Section:D4basis}

We now turn our attention to the dihedral symmetry group $D_4$ associated with a square. \addition{As mentioned in the introduction, w}hen we consider the generators \addition{of a $90^\circ$ degree rotation and a vertical reflection}
applied to the degree-$n$ monomial basis we induce a representation $\rho_n : D_4 \rightarrow O(n+1)$ with generators given by
\meeq{
\bfp_n(\addition{x,-y}) = \Vectt[x^n, -x^{n-1}y,\dots,(-y)^n] =  \underbrace{\diag(1, -1, \dots, (-1)^n)}_{\rho_n(s)}\bfp_n(\bfx) , \ccr
\bfp_n(\addition{-y,x}) = \Vectt[(-y)^n, x (-y)^{n-1},\dots,x^n] = \underbrace{\sopmatrix{ && (-1)^n \\ & \iddots \\ 1}}_{\rho_n(r)} \bfp_n(\bfx).
}


To deduce a symmetry-adapted basis we can either proceed numerically, or in this case we can deduce it explicitly. \addition{The final form of the basis is actually intuitive on its own, but in preparation for the cube where there is no intuitive explicit basis we present its construction in terms of explicit block-diagonalisation of the generators.}   The way to do this is probably best described by looking at some low degree examples. Constants are clearly invariant hence the representation is already irreducible, $\rho_0(g) = \addition{\sigmat}(g) = 1$, \addition{and we define $\bfq_0(x,y) := \bfp_0(x,y) = 1$ The} the linear case is precisely the irreducible faithful representation, $\rho_1(g) = \sigmaf(g) = g$ \addition{and we define $\bfq_1(x,y) := \bfp_1(x,y) = \vectt[x,y]$}.  For quadratics we compute
\[
Q_2 = \left( \begin{array} {c | c | c}
1/\sqrt{2} &  1/\sqrt{2} & \\
0  & 0 & 1\\
1/\sqrt{2} & - 1/\sqrt{2}
\end{array} \right)
\]
which diagonalises the representation, where the block structure is used to denote different irreducible representations.  In particular we have
\meeq{
Q_2^\top \rho_2(s) Q_2 = \diag(1,1,-1) = \diag(\addition{\sigmat}(s), \addition{\sigmats}(s), \addition{\sigmas}(s)) \ccr
Q_2^\top \rho_2(r) Q_2 = \diag(1,-1,-1) = \diag(\addition{\sigmat}(r), \addition{\sigmats}(r), \addition{\sigmas}(r)) 
}
Thus we have reduced the representation to a trivial, a reflect and a sign representation.  We deduce a symmetry-adapted basis as
\[
\bfq_2(x,y) := Q_2^\top  \bfp(x,y)  = \pr({\begin{array}{c|c|c}{{x^2 + y^2 \over \sqrt{2}}} & {x^2 - y^2  \over \sqrt{2}} & xy \end{array}})^\top.
\]
Continuing to the cubic case we can block-diagonalise via
\[
Q_3 = \left( \begin{array} {c c  c : c c  c}
1 & 0 &&& 0 & 0\\ 
0 & 0 &&& 0 & 1\\ 
0 & 0 &&& 1 & 0\\ 
0 & 1 &&& 0 & 0 
\end{array} \right)
\]
where the dash line is used to indicate that there are two copies of the same irreducible representation. 
We can deduce from checking the generators that 
\[
Q_3^\top \rho_3(g) Q_3 = \diag(\sigmaf(g) , \sigmaf(g)) = \sigmaf(g)^{\oplus 2},
\]
that is, we have two copies of a faithful representation. The corresponding basis is:
\[
\bfq_3(x,y) := Q_3^\top \bfp(x,y)  = \pr(\begin{array}{c c : c c }
x^3 & y^3 & x y^2 & x^2 y
\end{array})^\top.
\]
Unlike  \addition{bases corresponding to scalar irreducible representations these} work in pairs: applying a rotation does not leave $x^3$ invariant but rather leaves  ${\rm span}(x^3,y^3)$  invariant.

Our final example will be quartics, where we have all scalar irreducible representations present. In particular, we can block-diagonalise the representation using
%
%
\[
Q_4 = \pr( \begin{array} {c: c  | c  | c | c}
0 & 1/\sqrt{2} & 1/\sqrt{2} & 0 & 0 \\
0 & 0 & 0 &1/\sqrt{2} & 1/\sqrt{2} \\
1 & 0 & 0 & 0 & 0 \\
0 & 0 & 0 &-1/\sqrt{2} & 1/\sqrt{2} \\
0 &1/\sqrt{2} & -1/\sqrt{2} & 0 & 0 \\
\end{array}).
\]
This diagonalises the representation as
\[
Q_4^\top \rho_4(g) Q_4 = \diag(\addition{\sigmat}(g)^{\oplus 2}, \addition{\sigmats}(g), \addition{\sigmast}(g), \addition{\sigmas}(g)).
\]
and hence our symmetry-adapted basis is
\[
\bfq_4(x,y) :=\left( \begin{array}{c : c | c | c | c }
x^2 y^2 & {x^4+y^4 \over \sqrt{2}} & {x^4-y^4 \over \sqrt{2}}  & {x^3y - xy^3 \over \sqrt{2}} &  {x^3y + xy^3 \over \sqrt{2} }
\end{array}\right)^\top.
\]

\def\polyfourcoeffs{\left( \begin{array}{c : c | c | c | c }
3& {6 \over \sqrt{2}}&  -{4 \over \sqrt{2}}& -{2\over\sqrt{2}}& {6\over\sqrt{2}}
\end{array}\right)}
\def\polythreecoeffs{\pr(\begin{array}{c c : c c }
6 & 9 & 8  & 7
\end{array})}

As an example, consider the quartic  polynomial
\[
f(x,y) = x^4 + 2x^3y + 3x^2y^2 + 4xy^3 + 5y^4 \addition{+ 6x^3 + 7x^2y + 8 x y^2 + 9y^3}.
\]We  deduce that the coefficients in the symmetry-adapted basis are
\meeq{
\addition{\begin{pmatrix} 1 & 2 & 3& 4 & 5\end{pmatrix} Q_4}   =\addition{ \polyfourcoeffs}, \ccr
\addition{\begin{pmatrix} 6 & 7 & 8 & 9\end{pmatrix} Q_3}   =\addition{ \polythreecoeffs}.
}
Therefore we have
\meeq{
f(x,y) = \addition{\begin{pmatrix} 1 & 2 & 3& 4 & 5\end{pmatrix}  \bfp_4(x,y) + \begin{pmatrix} 6 & 7 & 8 & 9\end{pmatrix}  \bfp_3(x,y)} \ccr
 = \addition{\polyfourcoeffs \bfq_4(x,y) + \polythreecoeffs \bfq_3(x,y) } \ccr
= \ \underbrace{3x^2 y^2 +  3 (x^4 + y^4)}_{\hbox{Trivial}} \  + \  \underbrace{2 (y^4-x^4)}_{\hbox{Reflect}} \  + \  \underbrace{x y^3 - x^3 y}_{\hbox{Rotate}}\   +\   \underbrace{3 (x^3 y + xy^3)}_{\hbox{Sign}} \\
&\qquad \addition{+\  \underbrace{6 x^3 + 8 x y^2}_{\hbox{Faithful (1)}} \  + \   \underbrace{9 y^3 + 7 x^2 y}_{\hbox{Faithful (2)}}},
}
\addition{where there are two terms for the Faithful representation as it is two-dimensional. Note the symmetry of the faithful terms is hidden when we expand out the terms like this: it is only when left as coefficients in a symmetry-adapted basis that we can see that symmetry presents itself via the faithful representation.}

The pattern  for extending the construction is clear.  Even monomials can be combined to form irreducible symmetry-adapted bases corresponding to the scalar irreducible representations. For $\alpha = \floor{n/4}$ and $\beta = \floor{(n-2)/4}$ we have (up to multiplication by $1/\sqrt{2}$):
\[
\begin{array}{c | cc}
 \hbox{Trivial} & x^{n-2k}y^{2k} + x^{2k} y^{n-2k}, & k = 0,\ldots,\alpha,\\[3pt]
\hbox{Reflect}  &  x^{n-2k} y^{2k} - x^{2k} y^{n-2k}, & k = 0,\ldots,\beta,  \\[3pt]
\hbox{Rotate}  &  x^{n-1-2k} y^{2k+1} - x^{2k+1} y^{n-2k-1}, & k = 0,\ldots,\alpha-1,  \\[3pt]
\hbox{Sign} &  x^{n-1-2k} y^{2k+1} + x^{2k+1} y^{n-2k-1}, & k = 0,\ldots,\beta.  \\[3pt]
\end{array}
\]
Then for odd $n$ we have $(n+1)/2$ copies of the faithful representation with the basis
\[
\left( \begin{array}{c c} x^{n-2k} y^{2k} & x^{2k} y^{n-2k} \end{array} \right)^\top, \quad k =0,\ldots, (n-1)/2.
\]

\Subsection  Symmetry-adapted basis with respect to $O_h$.

We now turn our attention to constructing a symmetry-adapted basis in a cube with $O_h$ symmetries, again beginning with monomials which we write in lexographical order, i.e.,
\[
\bfp_n(x,y,z) = \vectt[x^n, x^{n-1} y, x^{n-1} z, \ldots, x y^{n-1}, x z^{n-1}, y^n, \ldots, y z^{n-1}, z^n].
\]
We will proceed by low degree examples but this can be generalised to arbitrary dimensions. 
The constant ($n = 0$) case is again the trivial representation and hence we have an irreducible symmetry-adapted basis $\bfq_0(x,y) := \bfp_0(x,y) = 1$. The $n = 1$ case gives us \addition{the generators}
\meeq{
\bfp_1(\tau_1 \bfx) = \vectt[z,-y,x] = \tau_1\bfp_1(\bfx)  \ccr
\bfp_1(\tau_2 \bfx) = \vectt[-x,-z,-y] = \tau_2\bfp_1(\bfx)  \ccr
\bfp_1(\tau_3 \bfx) = \vectt[-z,-y,-x] =  \tau_3 \bfp_1(\bfx)  \ccr
\bfp_1(s \bfx) = \vectt[-x,-y,-z] = s \bfp_1(\bfx),
}
i.e., we have $\rho_1(g) = g$.  This is an irreducible representation but not
in canonical form, which is necessary for \secref{PDEs}. Fortunately, an explicit orthogonal transformation to canonical form is presented in \secref{cubeirreducibles} which
we denote $Q_1$. Thus we have an irreducible symmetry-adapted basis $\bfq_1(x,y) := Q_1^\top \bfp_1(x,y)$.

For quadratics we have \addition{the generators}
\meeq{
\bfp_2(\tau_1 \bfx) = \vectt[z^2, -yz, xz, y^2,-yx,x^2] = \sopmatrix{ &&&&&1 \\ &&&&-1 \\ && 1 & 0 \\ &&0& 1 \\ &-1 \\ 1} \bfp_2(\bfx), \ccr
\bfp_2(\tau_2 \bfx) = \vectt[x^2, xz, xy, z^2, yz, y^2] = \sopmatrix{1 \\ &0& 1 \\ & 1 & 0 \\ &&&0&0& 1 \\ &&&0&1 & 0 \\ &&& 1 &0 &0 } \bfp_2(\bfx),  \ccr
\bfp_2(\tau_3 \bfx) =\vectt[z^2, yz, xz, y^2, xy, x^2] = \sopmatrix{ &&&&&1 \\ &&&&1 \\ && 1 & 0 \\ &&0& 1 \\ &1 \\ 1} \bfp_2(\bfx), \ccr
\bfp_2(s \bfx) = \bfp_2(\bfx).
}
This representation is still block-diagonalisable explicitly. In particular, if we define
\[
Q_2 := \pr(\begin{array}{cc | ccc|c}
	1/\sqrt{2} & 1/\sqrt{6} & 0 & 0 & 0 & -1/\sqrt{3}\\
	0 & 0 & 1/\sqrt{2} & -1/\sqrt{6} & 1/\sqrt{3} & 0\\
	0 & 0 & 0 & 2/\sqrt{6} & 1/\sqrt{3} & 0\\
	0 & -2/\sqrt{6} & 0 & 0 & 0 & -1/\sqrt{3}\\
	0 & 0 & 1/\sqrt{2} & 1/\sqrt{6} & -1/\sqrt{3} & 0\\
	-1/\sqrt{2} & 1/\sqrt{6} & 0 & 0 & 0 & -1/\sqrt{3}
\end{array})
\]
then we have
\[
Q_2^\top \rho_2(g)Q_2 = \diag(\sigma_{2+2,{\rm t}}(g), \sigma_{3+1,{\rm t}}(g), \sigma_{4,{\rm t}}(g))
\]
and we define $\bfq_2(x,y) := Q_2^\top \bfp_2(x,y)$. 

Going beyond quadratics we proceed numerically. In particular, we have $\bfp_n(s \bfx) = (-1)^n \bfp_n(\bfx)$  \addition{and multiplying the basis by an orthogonal matrix} will preserve this. 
We can deduce the generators $\rho_n(\tau_k)$  algorithmically which form a representation of $S_4$; in our implementation we do so by producing a vector of 3-tuples encoding the exponents of the monomials, apply the appropriate permutations arising from the change-of-variables $\tau_k$, and deduce the permutation that sorts this permuted vector of tuples. We can then deploy \cite{SOSymGroup} to compute a $Q_n$ (which is not uniquely defined!) that block-diagonalises the representation.
While it is not immediately clear which irreducible representations will be present they can be determined numerically, for example
\meeq{
Q_3^\top \rho_3(g) Q_3 = \diag(\sigma_{1+1+1+1,{\rm s}}(g), \sigma_{2+1+1,{\rm s}}(g)^{\oplus 2}, \sigma_{3+1,{\rm s}}(g)), \ccr
Q_4^\top \rho_4(g) Q_4 = \diag(\sigma_{2+1+1,{\rm t}}(g), \sigma_{2+2,{\rm t}}(g)^{\oplus 2}, \sigma_{3+1,{\rm t}}(g)^{\oplus 2}, \sigma_{4,{\rm t}}(g)^{\oplus 2}),\ccr
Q_5^\top \rho_5(g) Q_5 = \diag(\sigma_{1+1+1+1,{\rm s}}(g), \sigma_{2+1+1,{\rm s}}(g)^{\oplus 4},\sigma_{2+2,{\rm s}}(g), \sigma_{3+1,{\rm s}}(g)^{\oplus 2}), \ccr
Q_6^\top \rho_6(g) Q_6 = \diag(\sigma_{1+1+1+1,{\rm t}}(g),\sigma_{2+1+1,{\rm t}}(g)^{\oplus 2}, \sigma_{2+2,{\rm t}}(g)^{\oplus 3}, \sigma_{3+1,{\rm t}}(g)^{\oplus 4}, \sigma_{4,{\rm t}}(g)^{\oplus 3}).
}
We can keep going, and by degree $n=9$ we have \addition{a symmetry-adapted basis representing}  all 10 irreducible representations.

\Figurew[cubedecomposition]{\hsize}
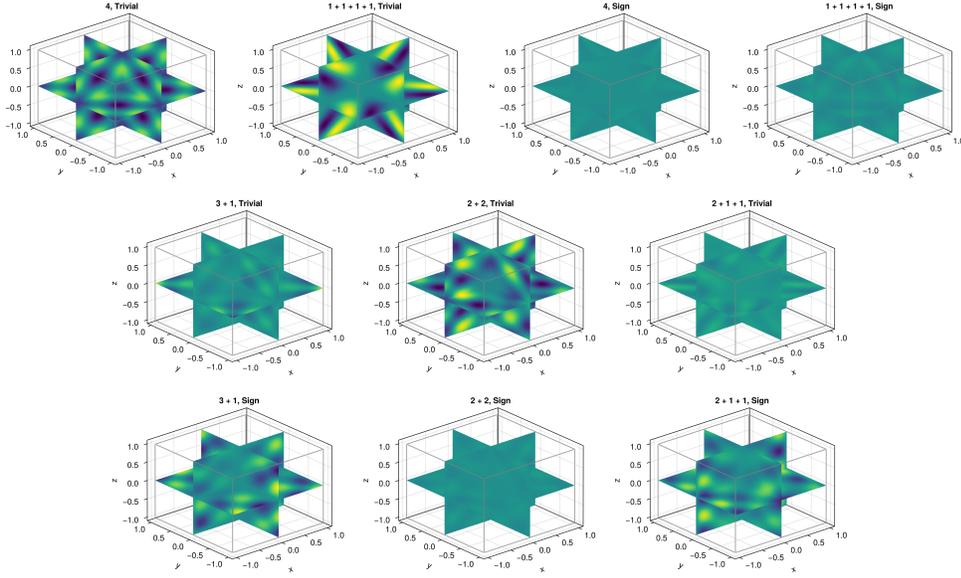
	A symmetric decomposition of the function $\sin(5(x-1/10)^2+y-3/10)\cos(6y+7z-2/5)$ on a cube into ten different symmetry classes, where for the matrix-valued representations we sum over the whole basis. The first row corresponds to scalar representations, the others are matrix representations and so the symmetry properties are not visually obvious but applying a symmetry action leaves each of these functions in their symmetry class.

In \figref{cubedecomposition} we consider decomposing an arbitrary  function into each symmetry class.  The first row correspond to scalar representations and any symmetry action leaves the function unchanged apart from possibly a sign change. The remaining 6 functions  correspond to matrix-valued irreducible representations where we sum over the whole basis; a more complete diagram would have a total of 20 different functions where each basis element of the irreducible representations are plotted separately.  \addition{Like the faithful terms in the square expansion, the symmetry is only apparent when expanded in a symmetry-adapted basis.}



\Section{PDEs}  Decoupling PDEs via symmetry-adapted bases.

%
%

We now return to the main theme of this paper: a  symmetry-adapted basis has no communication between bases corresponding to different irreducible representations when used to discretise PDEs. Moreover, because of Schur's lemma, there is no communication between different basis elements of each canonical irreducible representation. Thus we can split the numerical solution of a PDE in a square across 6 different independent linear systems (4 for each each scalar irreducible representation and another 2 for the two-dimensional faithful representation) and in a cube across 20 different systems.

Consider \addition{again} a Schr\"odinger equation
$
\Delta u + a(\bfx) u = \addition{f}
$
with \addition{Neumann} boundary conditions on a geometry $\Omega$ which is invariant under a symmetry action.
We focus on a weak-formulation of the equation though  everything we shall discuss extends to strong-formulation and more complicated PDEs provided some symmetry is present as outlined in the conclusions.  More precisely:  
 given finite polynomial vector spaces $U$ (the {\it trial space}) and $V$  (the {\it test space}), find $u \in U$ such that
\[
\ip<\nabla v, \nabla u> + \ip< v, a  u> =  \ip<v,\addition{f}>
\]
for all $v \in V$, where we assume the standard inner product
$\ip<f,g> := \int_\Omega f(\bfx) g(\bfx) \D \bfx$.

Given bases for $U$ and $V$ we arrive at a discretisation of the PDE. Monomials may be the immediate choice of basis though for reliable and stable computations it is preferable to use orthogonal polynomials, and in fact these can lead to sparse discretisations \cite{SOActa}, but the discussion that follows is independent of the choice of basis.
 In particular, if $U= {\rm span}(p_1,…,p_N)$ and  $V= {\rm span}(q_1,…,q_N)$,
we write these as a vector of functions $\bfp(\bfx) =\vectt[p_1(\bfx) , \dots , p_N(\bfx)]$
so that for $\bfc, \bfd \in \bbR^N$ we have
\meeq{
u(\bfx) = \sum_{k=1}^N c_k p_k(\bfx) = \bfp(\bfx)^\top \bfc, \quad v(\bfx) = \sum_{k=1}^N d_k q_k(\bfx) = \bfq(\bfx)^\top \bfd, \ccr
\addition{f(\bfx)} = \addition{\sum_{k=1}^N b_k p_k(\bfx) = \bfp(\bfx)^\top \bfb}.
}
Plugging this into the weak-formulation  reduces the problem to finding $\bfc \in \bbR^N$ such that
\[
\bfd^\top \pr( \ip<\nabla \bfq, \nabla \bfpt>  +  \ip< \bfq, a  \bfpt>) \bfc =   \bfd^\top \ip<\bfq,\bfpt> \addition{\bfb}
\]
for all $\bfd \in \bbR^N$, where here the inner products act entrywise, e.g.
\meeq{
 \ip<\bfq,a\bfpt>  = \sopmatrix{\ip<q_1,ap_1> & \cdots & \ip<q_1,ap_N>  \\ \vdots & \ddots & \vdots \\ \ip<q_N,ap_1> & \cdots & \ip<q_N,ap_N>}, \ccr
 \ip<\nabla\bfq,\nabla\bfpt>  = \sopmatrix{\ip<\nabla q_1,\nabla p_1> & \cdots & \ip<\nabla q_1,\nabla p_N>  \\ \vdots & \ddots & \vdots \\ \ip<\nabla q_N,\nabla p_1> & \cdots & \ip<\nabla q_N,\nabla p_N>}.
}
Since this holds true for all $\bfd$, including $\bfe_k$ for $k=1,\ldots,N$, we have arrived at \addition{a linear system}:
\[
\bigl(\underbrace{\ip<\nabla \bfq, \nabla \bfpt>}_{\hbox{stiffness matrix}}  \quad + \quad  \underbrace{ \ip< \bfq, a  \bfpt>}_{\hbox{mass matrix}} \bigr) \bfc =    \ip<\bfq,\bfpt> \addition{\bfb}.
\]
This discretisation will in general be dense (though it is sparse if we use a standard Finite Element Method basis or  more generally a $p$-Finite Element Method basis built from orthogonal polynomials \cite{babuska1981p}). 
However, if we build the bases using representation theory we can ensure sparsity and in fact decouple completely between different irreducible representations.
We now prove this for the stiffness and mass matrices.

\Lemma{intertwining}
Suppose $\bfq$ and $\bfp$ are  symmetry-\addition{respecting} bases with respect to representations $\rho_\addition{1}$ and $\rho_\addition{2}$, respectively, and $a : \Omega \rightarrow \bbR$ is invariant under the symmetry group $G \subset O(d)$. Then  $\ip< \bfq  , a \bfpt>$ is an intertwining matrix with respect to $\rho_\addition{1}$ and $\rho_\addition{2}$.

\Proof 

This result falls out of a couple observations. \addition{Recall} that for any orthogonal representation $\rho$ we have $\rho(g)^\top = \rho(g)^{-1} = \rho(g^{-1})$.
And for any orthogonal matrix $g \in O(d)$ that leaves  $\Omega$ invariant integration is   unchanged under an associated change-of-variables:
\[
\int_\Omega f(g \bfx) \D \bfx =  \int_\Omega f(\bfx) \D \bfx,
\]
which can be seen since the Jacobian of the change-of-variables is given by $|\det g| = 1$. 
We introduce the notation of applying a symmetry action $(g \star f)(\bfx) := f(g \bfx)$, so that a symmetry-adapted basis has the property, for $g \in G$,
\[
g\star \bfp = \rho_\addition{2}(g) \bfp\qand g \star \bfq = \rho_\addition{1}(g) \bfq
\]
and note that the invariance of integration for orthogonal change-of-variables implies that
\[
\ip<g\star  v, u> = \int_\Omega v(g \bfx) u(\bfx) \D \bfx  = \int_\Omega v( \bfx) u(g^\top \bfx) \D \bfx =   \ip<v, g^\top\!\star u>.
\]
Finally,  if $a$ is invariant under $G$ we have for any other function $f$ and for $g \in G$
\[
(g\star[a f])(\bfx) = a(g \bfx) f(g \bfx) = a(\bfx) f(g \bfx)  = a(\bfx) (g \star f)(\bfx).
\]
Putting everything together we have:
\meeq{
	\rho_\addition{1}(g)\ip< \bfq  , a \bfpt >   =  \ip< g\star \bfq  , a \bfpt >   =
\ip<  \bfq  , g^\top \!\star [a \bfpt] >\ccr
 = \ip<  \bfq  , a g^\top\! \star\bfpt >  =  \ip<  \bfq  , a   [\rho_\addition{2}(g^{-1})\bfp]^\top> =  \ip<  \bfq  , a   \bfpt > \rho_\addition{2}(g).
}
\mqed

%
%
%
%
%

The argument adapts to the stiffness matrix:

\Lemma{intertwiningstiff}
Suppose $\bfq$ and $\bfp$ are  symmetry-respecting bases with respect to $\rho_\addition{2}$ and $\rho_\addition{1}$, respectively. Then \addition{the stiffness matrix}  $\ip< \nabla\bfq  , \nabla \bfpt >$ is an intertwining matrix with respect to $\rho_\addition{1}$ and $\rho_\addition{2}$.

\Proof
In this proof it is convenient to use the convention $\nabla \bfp^\top : \Omega \rightarrow\bbR^{d \times n}$ arising from applying the gradient to each column of $\bfp^\top$ and we write $\nabla \bfq \equiv (\nabla \bfq^\top)^\top : \Omega \rightarrow \bbR^{n \times d}$. 
Note that the chain rule for gradients can be written as $\nabla[g \star u ](\bfx) = g^\top [g \star \nabla u](\bfx)$, or in particular
$\nabla[g\star \bfp^\top] =  g^\top (g\star\nabla \bfp^\top)$ and $\nabla[g\star \bfq] =  (g\star\nabla \bfq) g^\top$.
Finally since $\rho_\addition{1}(g) \bfq = g \star \bfq$ one can deduce that
\[
\rho_\addition{1}(g) \nabla \bfq = \nabla[g\star \bfq] = g^\top(g\star[\nabla \bfq]),
\]
with a similar formula relating $\rho_\addition{2}$ and $\nabla \bfp$.
Hence the result follows very similarly to the previous lemma:
\meeq{
\rho_\addition{1}(g) \ip<  \nabla\bfq  , \nabla \bfpt >  
 =   \ip<  g^\top (g\star\nabla\bfq)  , \nabla \bfpt >
 =   \ip<  \nabla\bfq  ,  g (g^\top\!\star \nabla \bfpt) >\ccr
 =   \ip<  \nabla\bfq  ,   \nabla[g^\top\!\star \bfpt] > 
 =   \ip<  \nabla\bfq  ,   \nabla \bfpt > \rho_\addition{2}(g).
}

\mqed

Thus for any symmetry-respecting basis our discretisation of a Schr\"odinger equation gives us intertwining matrices. But then Schur's lemma guarantees \addition{that} these decouple provided our basis is \addition{symmetry-adapted}!

\Theorem{blockdiag}
Suppose $\bfq$ and $\bfp$ are {\it symmetry-adapted bases}.  Then the mass and stiffness matrices are block-diagonalised after permutation, where  the number of blocks is at most the sum of the dimensions of all irreducible representations.

\Proof 
{
Enumerating the \addition{canonical} irreducible representations as $\sigma_1,\ldots,\sigma_r$ we can write the symmetry-adapted bases as
\meeq{
\bfq = \pr({\begin{array}{c:c:c|c|c:c:c}
(\bfq_1^{(1)})^\top& \ldots&(\bfq_1^{(a_1)})^\top &\ldots &(\bfq_r^{(1)})^\top&\ldots&   (\bfq_r^{(a_r)})^\top
\end{array}})^\top\ccr
\bfp = \pr({\begin{array}{c:c:c|c|c:c:c}
(\bfp_1^{(1)})^\top& \ldots&(\bfp_1^{(b_1)})^\top &\ldots &(\bfp_r^{(1)})^\top&\ldots&   (\bfp_r^{(b_r)})^\top
\end{array}})^\top
}
where 
$\bfq_\kappa^{(k)}$ and $\bfp_\ell^{(j)}$ are irreducible symmetry-adapted bases with respect to $\sigma_\kappa$ and $\sigma_\ell$, respectively.  \lmref{SchurMatrix} guarantees that if $\kappa \neq \ell$ then the associated mass matrix $\ip<\bfq_\kappa^{(k)}, a \bfp_\ell^{(j)}>$ is zero, as is the associated stiffness matrix. Thus the mass and stiffness matrix are block-diagonal with $r$ different blocks. But we also have $\ip<\bfq_\kappa^{(k)}, a \bfp_\kappa^{(j)}> = c_{\kappa,k,j} I_{d_\kappa}$ where $d_\kappa$ is the dimension of $\sigma_\kappa$ hence each of these blocks can be viewed as a block matrix where the blocks are constant multiples of the identity. Permuting the rows/columns to group together different basis elements results in a block-diagonal matrix with $d_\kappa$ blocks. Thus the total number of blocks is \addition{at most the sum of the dimensions} $d_1 + \cdots + d_r$. 
}
\mqed


\Figurew[squarespy]{\hsize}
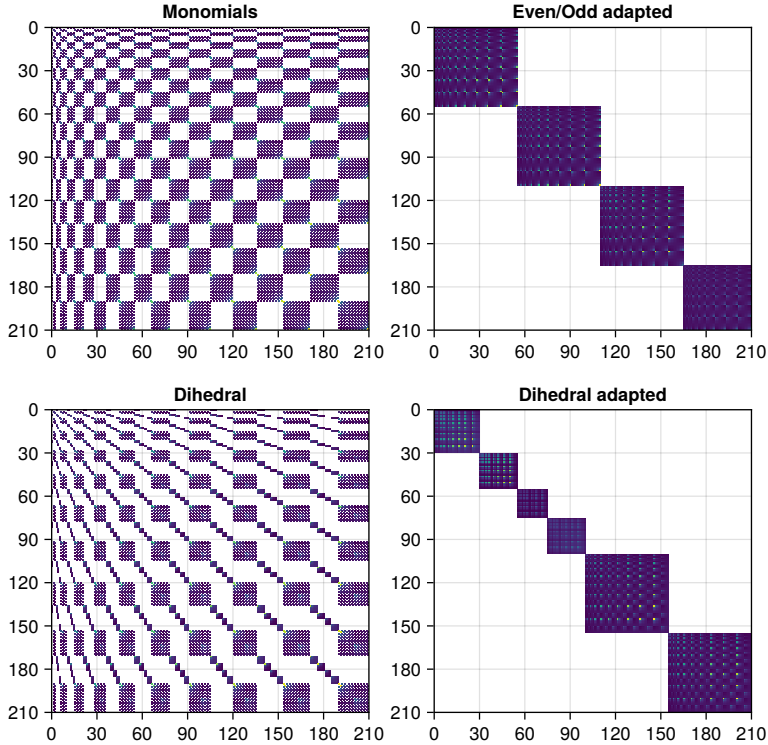
	Spy plots for the Schr\"odinger operator $\Delta u + (x^2 + y^2) u$ in the square discretised via weak-formulation \addition{using polynomials up to degree $19$}.  A monomial basis has a lot of sparsity in its discretisation (top left) and permuting rows/columns  according to even/odd symmetry  (that is, $\addition{\rmC}_2^2$) splits the discretisation into 4 independent matrices (top right). We have also construct an explicit basis respecting the $D_4$ symmetries resulting in increased sparsity (bottom left). The rows/columns can be  permuted to group together different basis elements of irreducible representations resulting in 6  independent matrices (bottom right). \addition{However, the dimensions of the blocks corresponding to the odd order polynomials are exactly the same for both the monomial and symmetry-adapted basis and hence there may be no efficiency gains when fully parallelised.}

\subsection{Examples}  First consider the Schr\"odinger equation in a square with potential $a(x,y) = x^2 + y^2$ which is invariant under the $D_4$ symmetry actions. In \figref{squarespy} we compare the sparsity pattern of the discretisation for the monomial basis which can be reordered according to even/odd ($\addition{\rmC}_2^2$) symmetry to result in 4 independent matrices, or with a symmetry-adapted basis for $D_4$ resulting in 6 independent matrices. This is a somewhat underwhelming example as the two sub-matrices associated with the faithful representations are just as big as those associated with $\addition{\rmC}_2^2$, so we would not expect any speed up if parallelised over 6 independent workers though may be faster if comparing parallelisation over 1, 2, or 3 workers. Note that the exact form of the potential does not actually change the sparsity pattern, provided it is invariant under $D_4$.

\Figurew[cubespy]{\hsize}
	A spy plot of the Schr\"odinger operator $\Delta u + (x^2 + y^2 + z^2) u$ in the cube discretised via weak formulation  \addition{using polynomials up to degree $9$}. The monomial basis splits according to even-odd symmetry (left) and a symmetry-adapted basis with respect to $O_h$ (right). We see an over double the number of decoupled matrices (from 8 to 20), with smaller sizes too which may lead to significantly faster simulations.

We now turn our attention to a 3D problem: consider the Schr\"odinger equation in a cube with potential $a(x,y,z) = x^2 + y^2 + z^2$ which is invariant under the $O_h$ symmetry actions. In \figref{cubespy} we  compare the sparsity pattern of using even-odd symmetry ($\addition{\rmC}_2^3$) which decouples across 8 different matrices with the full $O_h$ symmetry that decouples across 20 different matrices. In this case, not only do we get over double the amount of independent problems, each of the problems has notably smaller dimension than using only even-odd symmetry. \addition{In particular, when discretising with polynomials up to degree $9$, the largest block of the  permutation adapted discretisation is $22$ compared to $35$ for the monomial basis which only captures horizontal-vertical reflection, which means the discretisation size is reduced by about $63\%$.   For a dense solver that takes roughly $Cn^3$ operations one can therefore expect the solver to take as little as $\approx 0.63^3 \approx 1/4$ the time, that is, a $4\times$ speedup. Though this will vary in practice.}


\Figurew[cubepermspy]{\hsize}
	A spy plot of the Schr\"odinger operator $\Delta u + ((x-y)^2 + (y-z)^2 + (x-z)^2) u$ in the cube discretised via weak formulation \addition{with polynomials up to degree $9$} using the monomial basis splitting according to negation symmetry (left) and a symmetry-adapted basis with respect  to permutation and negation, i.e., $S_3 \times \addition{\rmC}_2$ (right).  We now see quadruple the number of decoupled matrices (from 2 to 8) which will lead to an increase in performance in simulations.  \addition{The worst case dimension of the blocks are now roughly a third the size leading to potentially as much as a 28$\times$ speedup. }
	
\Section{PermutationPDE} \addition{Schrödinger equations with multiple particles}.

There is another angle we can consider; rather than utilising more symmetry than even-odd symmetry we can consider potentials with less symmetry, e.g., consider the Schr\"odinger equation in a cube with potential $a(x,y,z) = (x-y)^2 + (y-z)^2 + (x-z)^2$ which is invariant under permuting $x,y,z$ and also invariant under negation $a(x,y,z) = a(-x,-y,-z)$. The solutions to the resulting equation can be interpreted as the wave function corresponding to three one-dimensional particles with quadratic interactions. Monomials capture the negation symmetry and enable decoupling into two linear systems. On the other hand, we can combine negation and permutation in the group $S_3 \times \addition{\rmC}_2$, very much as we did with $O_h \cong S_4 \times \addition{\rmC}_2$, but where the generators $\tau_1$ and $\tau_2$ correspond to the simple transpositions $(x,y,z) \mapsto (y,x,z)$ and $(x,y,z) \mapsto (x,z,y)$. Since the three irreducible representations of $S_3$ have a total dimension of 4, we get 8 distinct linear systems, \addition{as depicted in \figref{cubepermspy}}.  That is, we have quadrupled the number of independent matrices, each of a much smaller dimension!  \addition{Specifically, the largest block of the  permutation adapted discretisation is $41$ compared to $125$ for the monomial basis which only captures reflection symmetry, which means the discretisation is reduced by about $33\%$.   One can therefore expect a dense solver to take as little as $0.33^3 \approx 1/28$ the time, that is, a $28\times$ speedup. The computational savings actually improve at higher polynomial degrees, similar to the case of permutation symmetry which we will  observe  in \figref{numirreps}.}

\Figuretwo[numirreps,BosFermReduction]
	\addition{Left: the growth in the number of irreducible representations for $S_n$ (which is  equal to the number of partitions of $n$), and the sum of the dimensions of all irreducible representations of $S_n$. This shows that a 10-dimensional PDE with permutation and negation symmetry can in theory be parallelised across almost 20,000 cores. Right: The reduction in dimension of the independent blocks from a monomial basis  and a symmetry-adapted basis. Comparing the largest blocks we see an algebraic improvement with $n$. If we focus on specific symmetry classes such as Bosons (symmetric) and Fermions (anti-symmetric) this becomes a super-algebraic improvement}.

\addition{This potential increase in parallelisation will be even more pronounced in higher dimensional problems associated with more particles.  Consider an $n$-dimensional Schr\"odinger eigenvalue problem where the potential has permutation and negation symmetry, for example:
\[
\Delta u + \sum_{k<j} (x_k - x_j)^2 u = \lambda u
\]
where $\Delta = \sum_{k=1}^n \partial_{x_k}^2$  is the $n$-dimensional Laplacian, where the PDE is solved in the hypercube $[-1,1]^n$ and again for simplicity we assume  Neumann conditions.  As in the 3D case, $n$-dimensional monomials are a symmetry-respecting basis, and we can use \cite{SOSymGroup} to block-diagonalise the corresponding representation of $S_n$ and thereby construct a symmetry-adapted basis.    In this case, the number of irreducible representations  of $S_n$ equals the number of partitions of $n$ and the dimension of the irreducible representation is given by the hook length formula, see e.g. \cite{Okounkov}, but neither of these have a simple form. Nevertheless, they are computable directly and we see in \figref{numirreps} how the number of partitions grows super-algebraically whilst the sum of the their dimensions grows super-exponentially: thus there is a tremendous potential for parallelisation. The right-hand figure shows a comparison of the largest blocks of a permutation-adapted basis compared to monomials, which still only capture negation, showing a reduction in dimension that improves algebraically with $n$. In theory, for a PDE in 10-dimensions the discretisation using symmetry-adapted polynomials up to degree 6 are 0.5\% the size of the corresponding monomial discretisation, and so if we assume a dense eigenvalue algorithm takes roughly  $C n^3$ this would imply an 8 million times speedup! However, the dimensions of the discretisation are enormous, so even this theoretical speedup is not enough to overcome the curse of dimension. Thus it is unrealistic that it will scale much beyond 4D or 5D.
}

\Figurew[BosonFermion]{0.9\hsize}
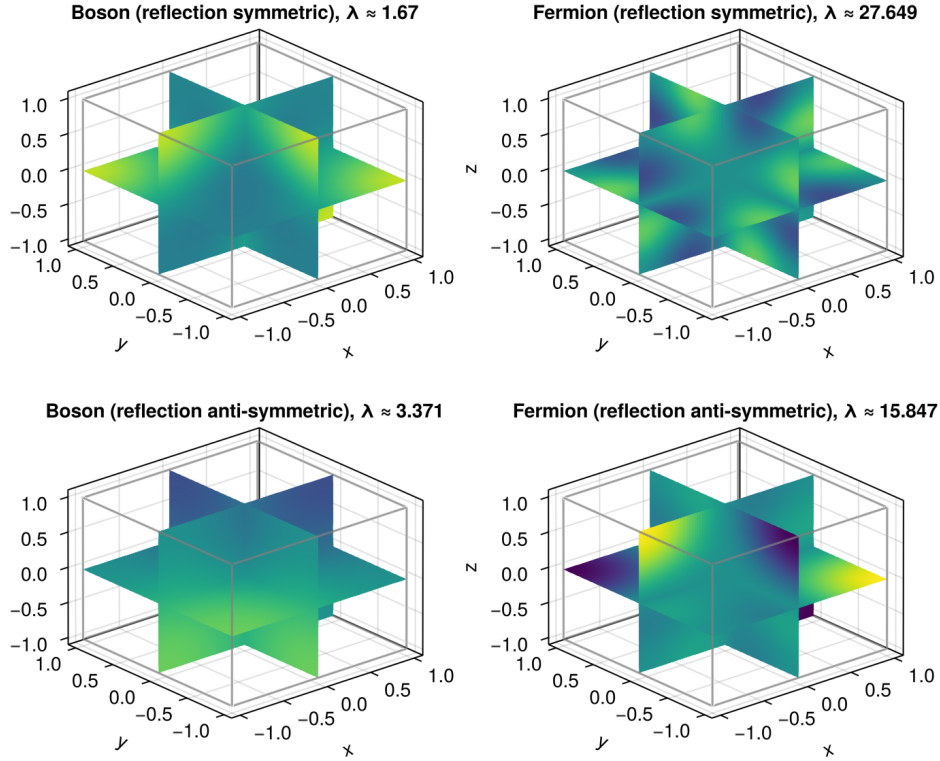
	Lowest energy Bosons and Fermions, where $\lambda$ is the associated eigenvalue of each of these eigenfunctions. Bosons are invariant under permutation of $x,y,z$ whilst Fermions are anti-invariant.  We have two types of Bosons and Fermions, each; those that are also invariant  under negation $(x,y,z) \mapsto (-x,-y,-z)$   and those that are anti-invariant.

One important feature in the quantum mechanical setting is that the type of solutions that are physically relevant are quite restricted: we want to compute Bosons which are invariant under permutation and Fermions\footnote{\addition{Fermions typically have a {\it spin}, which we have omitted as it is substantially more complicated. That is, we are considering  {\it spinless} Fermions.}} which are anti-invariant under permutation, which are directly related to the irreducible representations associated with trivial and sign representations. Because we also capture negation symmetry we can further decompose Bosons and Fermions into two distinct sub-classes. In \figref{BosonFermion} we depict four Bosons/Fermions corresponding to the minimal energy (eigenvalue) in each symmetry class.  By \addition{computing a symmetry-adapted basis corresponding to permutations and using only the terms corresponding to the trivial and sign representation  we can compute Bosons and Fermions directly, solving much smaller eigenvalue problems. For Fermions with three particles (i.e. in a cube) there is a reduction of dimension by around 5\% compared to a standard basis, which would lead to an expected $8000\times$ speedup for a dense solver.  This improvement is even more drastic in higher-dimensions: \figref{numirreps} shows the reduction in the size of the resulting systems, showing that they are at least exponentially smaller than a standard basis. This    may potentially overcome the curse of dimensionality (at least for one-dimensional particles), though the practical realisation of this discretisation requires more work: it is not realistic to build a very large representation corresponding to a tensor product basis and block-diagonalise it to determine the much smaller systems.  Alternatively,  Fermions can be numerically simulated using symmetry-adapted bases built from Slater derivatives in combination with sparse grids  \cite{griebel2007sparse} or directly in a spectral method \cite{clason2012general}, including in the important case of Coulomb interactions.
Very recently an approach based on symmetric polynomials has been used in Monte Carlo approximations in the many particle setting \cite{zhou2024multilevel}.
 }

\Section{conc} Conclusions.

Symmetry  can be used  to decouple \addition{discretisations} of partial differential equations into independent solves. This is standard practice for \addition{horizontal/vertical reflections (via even-odd bases), translations (via Fourier bases)} and rotations \addition{(via spherical harmonics)} but we have seen how \addition{discrete} symmetry groups like those associated with the square and cube can \addition{also} be used for  parallelisation. This is accomplished by \addition{incorporating} knowledge of irreducible representations \addition{in the construction of a symmetry-adapted basis, computed in practice via} either an explicit or numerical block-diagonalisation procedure. In the case of the cube \addition{with octohedral symmetry we can parallelise across 20 independent systems, or with permutation symmetry across 8 independent systems. When compared to using only reflection symmetry each linear system is notably smaller, hence can be solved significantly faster.} While we focussed on \addition{bases built from monomials} for simplicity everything translates to multivariate orthogonal polynomials, which on squares and cubes have the exact same \addition{construction procedure as the monomial basis.
While prior work in the 1970--2000s focussed on finite element methods where the techniques have limited practical implications, due in part to the success of multigrid, we view these techniques as a valuable tool in the arsenal of global spectral methods where multigrid is not available.
}
%
%

\addition{We explored the utilisation of these techniques for higher-dimensional Schr\"odinger equations where the potential for parallelisation and dimension reduction grows  with the dimension, particularly for the special cases of Bosons and (spinless) Fermions which have invariance and anti-invariance properties.} Note that approximating potentials invariant with respect to permutations by symmetric polynomials is an active research area with effective schemes recently introduced \cite{bachmayr2024polynomial}.  Such techniques have \addition{serious} potential to overcome the curse of dimensionality, \addition{though better techniques for constructing the resulting discretisation are needed to avoid prohibitively expensive pre-computation.}

A straightforward  extension is to vector-valued PDEs like Maxwell's equation, \addition{where the potential for efficiency gains in spectral methods is significant}. While we focussed on scalar \addition{PDEs} everything translates to the vector case if we modify the notion of a symmetry-adapted basis according to the principle of {\it equivariance}, that is, a   {\it  vector symmetry-respecting basis} $\bfp : \Omega \rightarrow \bbR^{n \times d}$ would satisfy
$
 \bfp(g \bfx) g  =   \rho(g) \bfp(\bfx)
$
and decomposing $\rho$ \addition{would produce} a vector symmetry-adapted basis. Weak formulations involving gradients, divergence,  and curl operators \addition{give rise to} intertwining matrices by very similar arguments to \lmref{intertwiningstiff}. Another example are Sturm--Liouville operators: if we have a matrix-valued function $A : \Omega \rightarrow \bbR^{d \times d}$ that is equivariant in the sense that, for all $g \in G$,
\[
A(g \bfx)  g = g A(\bfx)
\]
then \addition{discretising} an operator such as $\ip<\nabla v, A \nabla u>$ will also be an intertwining \addition{matrix}.
In all these cases Schur's lemma will guarantee that a discretisation automatically decouples according to the sum of the dimensions of the irreducible representations.

\addition{Finally, while at first glance these techniques appear limited to linear PDEs on highly specialised geometries and variable coefficients, it is possible that these simple problems may prove to be effective tools in preconditioners for more general PDEs. We point to recent work on Vertex Star Preconditioners by Brubeck and Farrell \cite{Brubeck2022} where parallelisable high-order numerical methods for simple (constant coefficient) PDEs are used to precondition more general problems. The proposed parallelisable techniques would widen the pool of preconditioners to allow variable coefficients with symmetry, though at this point the idea is merely speculative. }

\section*{Acknowledgments}
I thank Ryan Barnett, \addition{Marco Fasondini, Timon Gutleb}, Kaibo Hu, Martin Liebeck, Hans Munthe--Kaas, Peter Olver, \addition{Christoph Ortner}, Vic Reiner and \addition{Lior Silberman}  for helpful discussions. \addition{I also thank the anonymous referees for numerous suggestions which have substantially improved the paper.}

\bibliographystyle{siamplain}
\bibliography{RepTheoryPDEs}

\ends